\documentclass[11pt]{article}

\usepackage[a4paper,margin=1in]{geometry}
\usepackage{amsmath,amssymb,amsthm,mathtools}
\usepackage{enumitem}
\usepackage{aliascnt}
\usepackage{hyperref}
\usepackage[nameinlink,capitalize]{cleveref}
\hypersetup{colorlinks=true,linkcolor=blue,citecolor=blue,urlcolor=blue}

\newtheorem{theorem}{Theorem}[section]

\newaliascnt{proposition}{theorem}
\newtheorem{proposition}[proposition]{Proposition}
\aliascntresetthe{proposition}
\newaliascnt{lemma}{theorem}
\newtheorem{lemma}[lemma]{Lemma}
\aliascntresetthe{lemma}
\newaliascnt{corollary}{theorem}
\newtheorem{corollary}[corollary]{Corollary}
\aliascntresetthe{corollary}
\newaliascnt{problem}{theorem}
\newtheorem{problem}[problem]{Problem}
\aliascntresetthe{problem}

\theoremstyle{definition}
\newaliascnt{definition}{theorem}

\aliascntresetthe{definition}

\theoremstyle{remark}
\newaliascnt{remark}{theorem}
\newtheorem{remark}[remark]{Remark}
\aliascntresetthe{remark}
\newaliascnt{example}{theorem}
\newtheorem{example}[example]{Example}
\aliascntresetthe{example}

\crefname{theorem}{Theorem}{Theorems}
\crefname{proposition}{Proposition}{Propositions}
\crefname{lemma}{Lemma}{Lemmas}
\crefname{corollary}{Corollary}{Corollaries}
\crefname{problem}{Problem}{Problems}
\crefname{definition}{Definition}{Definitions}
\crefname{remark}{Remark}{Remarks}
\crefname{example}{Example}{Examples}

\newcommand{\End}{\operatorname{End}}
\newcommand{\Endo}{\operatorname{End}^{0}}
\newcommand{\Cent}{\operatorname{Cent}}
\newcommand{\Gal}{\operatorname{Gal}}
\newcommand{\Frac}{\operatorname{Frac}}
\newcommand{\coker}{\operatorname{coker}}
\newcommand{\rk}{\operatorname{rank}}
\newcommand{\cO}{\mathcal O}
\newcommand{\cQ}{\mathcal Q}
\newcommand{\Sfr}{\mathfrak S}
\newcommand{\pfr}{\mathfrak p}
\newcommand{\Q}{\mathbb Q}
\newcommand{\Z}{\mathbb Z}
\newcommand{\F}{\mathbb F}
\newcommand{\C}{\mathbb C}
\newcommand{\inv}{\operatorname{inv}}
\newcommand{\Spec}{\operatorname{Spec}}

\title{Simple abelian varieties with prescribed endomorphism ratio}
\author{Saeed Tafazolian}
\date{}

\hypersetup{pdftitle={Simple abelian varieties with prescribed endomorphism ratio},
  pdfauthor={Saeed Tafazolian},
  pdfkeywords={Abelian varieties, endomorphism algebras, p-rank zero, Honda-Tate theory}}

\begin{document}
\maketitle

\begin{abstract}
Let $A$ be a simple abelian variety over an algebraically closed field, and set
$m(A)=2\dim A/[\operatorname{End}^0(A):\mathbb Q]$.
We prove that, for every prime $p$ and every positive rational number $\mu$, there exist simple abelian varieties of characteristic $p$ with endomorphism ratio $\mu$, both of $p$-rank zero and of positive $p$-rank. In particular, we answer a question of Oort in a stronger form.

Our construction combines Honda–Tate theory with explicit deformations of $p$-divisible groups using Breuil windows. A Frobenius matrix with quaternionic multiplication gives two reduced generic slopes with coprime denominators. Specialization and constancy of dimension then determine the full geometric endomorphism algebra. Polarized Serre–Tate theory and Grothendieck algebraization lift these deformations to abelian schemes. A related construction with an additional étale factor yields the positive-$p$-rank result.
\end{abstract}



\section{Introduction}
For an abelian variety $A/k$, write
\[
\Endo_k(A)=\End_k(A)\otimes_{\mathbb Z}\mathbb Q,
\qquad
m_k(A)=\frac{2\dim A}{[\Endo_k(A):\mathbb Q]}.
\]
We omit the subscript when the base field is algebraically closed; thus all endomorphism algebras in the main theorem are geometric. In characteristic $p$, the $p$-rank $f(A)$ is defined by
\[
\#A[p](\bar{k})=p^{f(A)}.
\]

If $A$ is simple in characteristic zero, then $m(A)$ is a positive integer, and every positive integer occurs \cite[\S22]{OortQ}. In positive characteristic, however, $m(A)$ need not be integral. The numerical restrictions on endomorphism algebras \cite[\S21]{Mumford} imply that $m(A)$ is an integer for Albert types~I and~II, while for type~III its reduced denominator divides $2$. For type~IV, let $D=\Endo(A)$ and write
\[
[Z(D):\mathbb Q]=2e_0,
\qquad
[D:Z(D)]=d^2.
\]
Since $e_0d\mid\dim A$, one has $m(A)\in d^{-1}\mathbb Z$. The possible endomorphism algebras in positive characteristic are described in \cite{OortEnd1988,OortVanDerPut}.

Oort asked which values occur.  The question goes back to his 1995
list \cite[Question~2]{OortSomeQ}; it is formulated as an expectation,
for every $\gamma\in\Q_{>0}$ and every prime $p$, in
\cite[22.8--22.9]{OortFF}, and it appears in Appendix~3 of
\emph{Open Problems in Arithmetic Algebraic Geometry} in the following
form \cite[Question~22.2, p.~315]{OortQ}.

\begin{problem}\label{prob:oort-main}
Suppose given a positive rational number $\mu\in\Q_{>0}$ and a prime
number $p$.  Does there exist a simple abelian variety $A$ over an
algebraically closed field $k\supset\F_p$ such that $m(A)=\mu$?
\end{problem}

Oort adds: ``I expect every positive rational number does appear in
this way in every positive characteristic''
\cite[p.~315]{OortQ}.  We confirm this expectation in the following
stronger form.

\begin{theorem}\label{thm:main}
Let $p$ be a prime and $\mu\in\Q_{>0}$.  There exist an algebraically
closed field $k$ of characteristic $p$ and a simple abelian variety
$A/k$ such that $f(A)=0$ and $m(A)=\mu$.
\end{theorem}

The field $k$ can be taken to be an algebraic closure of a finitely
generated extension of $\F_p$ (\cref{rem:field-of-definition}).  The
construction is effective: for $\mu=a/b$ with integers $a\ge2$ and
$b\ge1$, it produces a variety of dimension $32ab$, and
\cref{cor:ratio-dimension} gives smaller dimensions for many ratios.
These are the dimensions produced by the construction; we make no claim
about optimality (\cref{rem:dimension}).  In \cref{sec:complements} we
show that every $\mu\in\Q_{>0}$ is also realized by a simple abelian
variety of positive $p$-rank (\cref{cor:both-ranks}).

For a simple abelian variety over an algebraic extension of $\F_p$,
descent to a finite field and the Honda--Tate dimension formula give
$m(A)=1/d$, where $d$ is the degree of its endomorphism division algebra
over its centre \cite[\S3, Theorem~2 and (7)]{Tate66};
see also \cite[Part~I, Theorem~8]{WaterhouseMilne}.  Every unit
fraction occurs: ordinary and supersingular elliptic curves give $1$
and $1/2$, respectively, and the quadratic seeds of \cref{prop:seeds}
give $1/d$ for $d\ge3$.  Consequently, if $\mu=a/b$ is in lowest terms
and $a>1$, no simple abelian variety realizing $\mu$ can be defined
over an algebraic extension of $\F_p$.

For a polarized abelian variety $Z/k$ and a Rosati-stable
$\Q$-subalgebra $D\subseteq\Endo_k(Z)$, Theorem~1.1 of
\cite{OortVanDerPut} constructs, for each $M\ge1$, an abelian variety
$Y/k((t))$ with geometric endomorphism algebra $D$ and
$\dim Y=\dim Z+M[D:\Q]$.  We record the consequence used for
comparison.

\begin{corollary}\label{cor:known-large-ratios}
For every prime $p$ and every rational $\mu>2$, there exists a simple
abelian variety $A$ over an algebraically closed field of
characteristic $p$ with $m(A)=\mu$.
\end{corollary}

\begin{proof}
Choose $b\ge3$ such that $N=b(\mu-2)$ is a positive integer.  The seed
$B=B_{p,b}$ of \cref{prop:seeds} has dimension $b$ and endomorphism
division algebra $D=D_{p,b}$ of $\Q$-dimension $2b^2$.  Apply
\cite[Theorem~1.1]{OortVanDerPut} to $Z=B^N$, its product polarization,
the diagonal copy of $D$, and $M=1$.  The geometric endomorphism algebra
of the resulting $Y$ is the division algebra $D$, so its geometric fibre
is simple, and
\[
             m(Y_{\overline{k((t))}})
             =\frac{2(Nb+2b^2)}{2b^2}=2+\frac Nb=\mu.
\]
\end{proof}

The ratio $2$ is realized by an elliptic curve with transcendental
$j$-invariant.  The ratios not realized by these classical constructions
are therefore the fractions $a/b$ in lowest terms with $b\ge2$ and
$2\le a<2b$.  Moreover, the Oort--van der Put construction produces
varieties of positive $p$-rank: part~(4) of their theorem identifies the
formal group of $Y$ with that of $Z$ times
$\widehat{\mathbb G}_m^{\,M[D:\Q]}$, and since the $p$-divisible group
of a polarized abelian variety is self-dual up to isogeny,
$f(Y)=f(Z)+M[D:\Q]>0$.  With $p$-rank zero, the classical constructions
recalled in \cref{rem:classical-prank-zero} realize the unit fractions,
the ratio $1$, and the even integers $2g\ge6$.

For the remaining ratios write $\mu=a/b$ with $a\ge2$; the same construction also gives $\mu=1$ by taking $a=b=2$. We realize $\mu$ as $n/d$, where $d\ge4$ is even and $n=r+s$ with coprime $s>r\ge3$. The condition that $d$ be even ensures that the prescribed local invariants have common denominator $d$, while $d\ge4$ is needed for the quartic CM seed constructed in Section~\ref{sec:seeds}.
We fix a quartic CM field $K$ in which $p$ splits completely and a
geometrically simple Honda--Tate variety $B_d$ with endomorphism
division algebra $D/K$ of degree $d$.  At one conjugate pair of primes
above $p$ the algebra $D$ has local division index $2$, and this pair
carries a quaternionic deformation; at the other pair the local index is
$d$, which fixes the global degree.  An explicit Frobenius matrix
gives a deformation of the first pair whose geometric generic fibre has
two reduced slopes, with denominators $r$ and $s$.  The matrix
interpolates linearly between two monomial matrices: one gives the
isoclinic special fibre, the other already has the desired generic
Newton polygon.  With polarized Serre--Tate theory we globalize the
construction to an abelian scheme with special fibre isogenous to
$B_d^n$.

To determine the geometric generic endomorphism algebra $E$, consider
the commutant $C=\Cent_E(D)$ of $D$ in $E$.  A completion of $C$ embeds
into a product of two division algebras of $\Q_p$-dimensions $r^2$ and
$s^2$.  A nontrivial idempotent would, after finite base change and
isogeny, give a summand with different dimensions in the special and
generic fibres, contradicting constancy of dimension.  The completed
commutant is therefore a division algebra, whose dimension divides both
$r^2$ and $s^2$.  Since $\gcd(r,s)=1$, it equals $\Q_p$, and the
double-centralizer theorem gives $E=D$ (\cref{prop:commutant-criterion}).
The resulting variety $A$ has $\dim A=2nd$ and $[D:\Q]=4d^2$, so
$m(A)=n/d=a/b$, and all its slopes lie strictly between zero and one.

Two parts of the argument may be useful elsewhere.  The commutant
criterion determines a geometric endomorphism algebra exactly from a
change of Newton polygon at a single prime above $p$; its local
hypotheses concern a single local factor, and none of its hypotheses
depends on the way the deformation is constructed.  The windows of Section~\ref{sec:explicit-matrix}, obtained
by interpolating between two monomial Frobenius matrices, carry an
action of a maximal order and have generic Newton polygons that can be
read off from a graph, uniformly in $p$, including $p=2$.  Combined with
CM fields in which $p$ splits completely, they let one pair of primes
carry the deformation while other pairs fix the global degree or, in
Section~\ref{sec:complements}, contribute an \'etale part.

Question~22.3 of \cite[p.~315]{OortQ} asks for a more specific
deformation.  In our notation (the source writes $a/d$, $A$, $B=A^a$
and $D:=\End(A)$), it reads as follows.  Fix a reduced fraction $a/b$
with $b>1$ and an abelian variety $B$ with $m(B)=1/b$, and equip $B^a$
with the diagonal polarization $\lambda$ and the diagonal action
$\iota:\End(B)\hookrightarrow\End(B^a)$.  Does the triple
$(B^a,\lambda,\iota)$ admit a deformation whose generic fibre has
endomorphism algebra isomorphic to $\End(B)$?  We read the last
condition as asking that the geometric generic fibre have endomorphism
algebra isomorphic to $\Endo(B)$.  Our construction uses powers $B_d^n$
with $n/d=a/b$ and allows a $p$-power isogeny of the special fibre, a
transported and scaled polarization, and an action of a chosen full
order in $\Endo(B_d)$.  It therefore does not give an affirmative answer
to Question~22.3 in the form stated there, with its prescribed special
fibre, polarization and action of the full endomorphism ring.  The
ranges in which our special fibre is isogenous to the corresponding
power of a unit-fraction seed, subject to these modifications, are
described in Section~\ref{sec:realization}.

Section~\ref{sec:strategy} establishes the commutant criterion.
Section~\ref{sec:explicit-matrix} constructs the Frobenius matrix and its
two integral realizations.  Section~\ref{sec:seeds} gives the
Honda--Tate seeds, and Section~\ref{sec:globalization} globalizes the
local deformations.  Section~\ref{sec:realization} proves
\cref{thm:main}, records the choice of parameters, and gives a smaller
quadratic-CM construction.  Section~\ref{sec:complements} realizes every
ratio with positive $p$-rank.

\section*{AI Declaration}

ChatGPT (OpenAI, GPT-6) was used to assist with literature searches, language revision, and the checking of intermediate computations and exposition. It occasionally suggested alternative formulations or simplifications. All mathematical arguments, calculations, references, and final formulations were independently checked and revised by the author, who takes full responsibility for the contents of the manuscript.

\section{The commutant criterion}\label{sec:strategy}

The key step excludes extra endomorphisms by comparing the dimensions
of slope summands on two fibres.  We state it independently of the local
matrix construction.

For a $p$-divisible group $X$, write
$\Endo(X)=\End(X)\otimes_{\Z_p}\Q_p$.  When the subscript on $\Endo$ is
an algebra acting on $X$, rather than a field of definition as in the
introduction, it specifies the algebra whose action is to be
centralized.  We use covariant
Dieudonn\'e theory: the Newton slopes of a $p$-divisible group lie in $[0,1]$ and their sum, counted with multiplicity, is its height minus its dimension, and Cartier duality replaces each slope
$\lambda$ by $1-\lambda$.  Thus an isoclinic group of slope $\lambda$
and height $\beta$ has dimension $\beta(1-\lambda)$.  Height and dimension are locally constant
in families and invariant under isogeny.  For an abelian variety $A$, the $p$-rank is the multiplicity of slope $1$ in the covariant Newton polygon of $A[p^\infty]$.

\begin{lemma}\label{lem:spread-idempotent}
Let $k$ be a perfect field of characteristic $p$, put $R=k[[y]]$, and
let $X/R$ be a $p$-divisible group with generic fibre $X_\eta$ over
$\Frac R$.  Every idempotent of $\Endo(X_\eta)$ becomes an integral
idempotent over $R$ after replacing $X$ by an isogenous $p$-divisible
group.
\end{lemma}

\begin{proof}
For an idempotent $e$, choose $N$ so that $p^Ne$ and $p^N(1-e)$ are
integral.  Zink's small-image construction
\cite[Proposition~8]{ZinkSlope} gives their images $Y_\eta$ and
$Z_\eta$ over $\Frac R$ and a map $v:X_\eta\to Y_\eta\oplus Z_\eta$.
The construction commutes with field extension.  After extension to an
algebraic closure, the rational Dieudonn\'e module of $v$ is
$p^N(e,1-e)$ onto complementary summands.  Hence $v$ is an isogeny,
also over $\Frac R$ by faithfully flat descent.

Its kernel is contained in $X_\eta[p^m]$ for some $m$.  The schematic
closure $H$ in $X[p^m]$ is finite flat over $R$. Indeed, its coordinate
ring is finite over $R$ and injects into the coordinate algebra of the
generic fibre, hence is torsion-free over the DVR.  It is a subgroup: the
group operations preserve it on the schematically dense generic fibre,
and the same holds for $H\times_RH$.  The quotient $X'=X/H$ is a
$p$-divisible group with $X'_\eta\simeq Y_\eta\oplus Z_\eta$.  The
transported projector is therefore an integral endomorphism of
$X'_\eta$.

Since $R=k[[y]]$ is a discrete valuation ring of characteristic $p$,
restriction of $p$-divisible groups from $R$ to its fraction field is
fully faithful by \cite[Corollary~1.2]{deJongHom}.  The projector extends uniquely to
$X'$, and its extension is idempotent by uniqueness.  An integral
idempotent splits a $p$-divisible group into its image and complementary
image, as required.
\end{proof}

We use standard consequences of Poincar\'e complete reducibility and
central simple algebra theory: geometric endomorphism algebras are
semisimple, centralizers of separable subalgebras in semisimple algebras
are semisimple, and Skolem--Noether and the double-centralizer theorem
apply to central simple algebras over a common centre; see
\cite[\S19]{Mumford}, \cite[\S2.1 and Theorem~2.7.2]{GilleSzamuely}
for Wedderburn's theorem and Skolem--Noether, and \cite[\S12.7]{Pierce}
for the double-centralizer theorem.
These facts will be used without further proof.

\begin{proposition}\label{prop:commutant-criterion}
Let $k$ be algebraically closed of characteristic $p$, put $R=k[[t]]$,
and let $\mathcal A/R$ be an abelian scheme with a rational action of
a central division algebra $D$ over a number field $K$, that is, a
unital homomorphism $D\to\Endo(\mathcal A)$.  Suppose that
the special fibre is $D$-linearly isogenous to $B^n$, with
$\Endo(B)=D$ and diagonal $D$-action.  Let $\mathfrak p\mid p$ satisfy
$K_{\mathfrak p}=\Q_p$, put $D_{\mathfrak p}=D\otimes_KK_{\mathfrak p}$,
and suppose that the corresponding rational central factor of
$\mathcal A[p^\infty]$ is represented by a direct factor $Y/R$.

Assume that $Y_0$ is isoclinic of slope $\lambda_0$, whereas
$Y_{\bar\eta}$ is isogenous to $Y_1\oplus Y_2$, where $Y_i$ is
isoclinic of slope $\lambda_i\ne\lambda_0$.  Assume further that,
under this decomposition,
\[
 \Endo_{D_{\mathfrak p}}(Y_{\bar\eta})
   \simeq \mathcal D_r\times\mathcal D_s,
 \qquad \dim_{\Q_p}\mathcal D_r=r^2,\quad
 \dim_{\Q_p}\mathcal D_s=s^2,
\]
where the two factors are division algebras acting on $Y_1$ and $Y_2$,
respectively, and $\gcd(r,s)=1$.  Then
$\Endo(\mathcal A_{\bar\eta})=D$, and the geometric generic fibre is
simple.
\end{proposition}

\begin{proof}
\smallskip\noindent\emph{Step 1. Specialization and the completed commutant.}
Put $E=\Endo(\mathcal A_{\bar\eta})$.  All elements of a $\Q$-basis of
$E$ are defined over a finite extension of $k((t))$.  Homomorphisms of abelian schemes extend
uniquely over normal bases \cite[Proposition~I.2.7]{FaltingsChai}, and
specialization is injective by rigidity.  We therefore obtain
\[
               D\subseteq E\hookrightarrow M_n(D).
\]
The centre $K$ specializes to scalar matrices; injectivity applied to
commutators shows that $K\subseteq Z(E)$.  Consequently
$C=\Cent_E(D)$ is a semisimple $K$-algebra and specializes into
$\Cent_{M_n(D)}(D)=M_n(K)$.

Set $\widehat C=C\otimes_KK_{\mathfrak p}$.  We show that the
completed commutant still acts faithfully.  Choose a $K$-basis
$c_1,\ldots,c_e$ of $C$, defined over one finite extension
$\Omega'/k((t))$, and let $R'$ be its complete valuation ring.  Its
residue field is $k$, and Cohen's structure theorem gives
$R'\simeq k[[y]]$, so \cref{lem:spread-idempotent} applies over $R'$.  After clearing a common
denominator, the $c_i$ extend to $\mathcal A_{R'}$ and preserve its
rational central factors.  On $Y_{R'}$, the field $K$ acts through
$K\to K_{\mathfrak p}=\Q_p$.  The restrictions of the $c_i$ to
$Y_{R'}$ therefore give
\[
 \rho:\widehat C\longrightarrow\Endo(Y_{R'}),\qquad
 \sum_i c_i\otimes a_i\longmapsto\sum_i a_i(c_i|_{Y_{R'}}).
\]
Here every $a_i\in\Q_p$ acts over the same base through the canonical
$\Z_p$-action of a $p$-divisible group.  Thus all elements of $\widehat C$,
including idempotents that do not descend to $C$, are defined over this
single extension.

On the special fibre, $\rho$ is the composite of the injection
$\widehat C\hookrightarrow M_n(K_{\mathfrak p})$ and the faithful
multiplicity-space action on the active factor of $B^n[p^\infty]$.
It is therefore injective.  Restriction to the generic fibre is
injective on homomorphisms of finite flat group schemes, and hence on
quasi-endomorphisms of $p$-divisible groups.  We obtain a unital embedding
\[
 \widehat C\hookrightarrow\Endo_{D_{\mathfrak p}}(Y_{\bar\eta})
     \simeq\mathcal D_r\times\mathcal D_s.
\]

\smallskip\noindent\emph{Step 2. Constancy of dimension excludes projectors.}
A nontrivial idempotent of $\widehat C$ maps to $(1,0)$ or $(0,1)$, because both
factors are division algebras.  By \cref{lem:spread-idempotent}, after
an isogeny over $R'$ it defines a direct summand of height
$\beta_i>0$ whose generic slope is $\lambda_i$.  Its special slope is
$\lambda_0$, since the special fibre of $Y$ is isoclinic.  Its
dimensions would therefore be
\[
                  \beta_i(1-\lambda_i)
                  \quad\text{and}\quad \beta_i(1-\lambda_0),
\]
which are different.  This contradicts constancy of dimension over the
connected base $\Spec R'$.  Thus $\widehat C$ has no nontrivial idempotents.
It is semisimple, since it is a characteristic-zero scalar extension of
the semisimple algebra $C$, and hence it is a division algebra.

\smallskip\noindent\emph{Step 3. Coprimality determines the endomorphism algebra.}
Both coordinate projections of $\widehat C$ are unital and therefore injective.
Viewing each target as a right vector space over $\widehat C$ gives
$[\widehat C:\Q_p]\mid r^2$ and $[\widehat C:\Q_p]\mid s^2$.  Coprimality implies
$[C:K]=[\widehat C:\Q_p]=1$, so $C=K$.
Every central idempotent of $E$ lies in $C$; thus $E$ is simple, and
$K\subseteq Z(E)\subseteq C=K$ makes it central simple over $K$.
The double-centralizer formula gives
$[E:K]=[D:K][C:K]=[D:K]$.  Hence $E=D$, which also proves simplicity
of $\mathcal A_{\bar\eta}$.
\end{proof}

In the proof of \cref{thm:main}, we use special slope $1/2$ and generic
slopes $1/2-1/r$ and $1/2+1/s$.  For an integer $a\ge2$, the choice
$(r,s)=(2a-1,2a+1)$ gives coprime centralizer degrees and keeps
these slopes in $(0,1)$.

\section{Explicit Frobenius operators and windows}\label{sec:explicit-matrix}

We now construct the local deformations needed for
\cref{prop:commutant-criterion}.  The same explicit Frobenius matrix will
be realized on two lattice chains.  We begin by explaining how to read
the geometric generic isocrystal from a window.

\subsection{Windows and generic slopes}\label{sec:deformation-prerequisites}

Let \(k\) be a perfect field of characteristic \(p\), put
\(\Sfr=W(k)[[t]]\) and \(R=\Sfr/(p)=k[[t]]\), and extend Witt
Frobenius by \(\sigma(t)=t^p\).  A Breuil window is a finite free
\(\Sfr\)-module \(Q\) with an injective map
\(\phi:Q\to Q^{(\sigma)}\) whose cokernel is annihilated by \(p\)
and is finite free over \(R\).  For the frame
\(\mathcal B=(\Sfr,p\Sfr,R,\sigma,\sigma_1)\), where
\(\sigma_1(px)=\sigma(x)\), we are in the setting of
\cite[\S6]{Lau2014} with \(E=p\), a power series with constant term
\(p\), and \(R=\Sfr/E\Sfr\).  Since \(\sigma(t)=t^p\), the conditions of
\cite[Proposition~6.2]{Lau2014} hold by \cite[Remark~6.3]{Lau2014}, so
\cite[Theorem~6.6 and Corollary~6.7]{Lau2014} give an equivalence
between \(p\)-divisible groups over \(R\) and Breuil windows, for every
prime \(p\), including \(p=2\).
The associated \(p\)-divisible group has height \(\rk Q\), dimension
\(\rk_R\coker(\phi)\), and covariant rational Frobenius with
linearization \(F^\sharp=p\phi^{-1}\).  In matrix formulas we write
\(F=p\phi^{-1}\), with the Frobenius twist understood.
Displays over $R$ are taken in the projective-limit sense of
\cite[\S2G, Lemma~2.16]{Lau2014}; we distinguish Zink's
ring $\mathbb W(R)$ from the ordinary Witt ring $W(R)$.

Let \(\varkappa\colon\Sfr\to\mathbb W(R)\subset W(R)\) be Lau's
coefficient map; thus \(\varkappa(t)=[t]\) and
\(F\varkappa=\varkappa\sigma\).  Put \(\Omega=\overline{k((t))}\).
Composing \(\varkappa\) with the Witt functoriality of \(R\to k\),
\(t\mapsto0\), and of the inclusion \(R\hookrightarrow\Omega\) gives
ring homomorphisms \(\kappa_0\colon\Sfr\to W(k)\) and
\(\kappa_\eta\colon\Sfr\to W(\Omega)\), with \(\kappa_0(t)=0\) and
\(\kappa_\eta(t)=[t]\), both satisfying \(\kappa\sigma=F\kappa\).
Applied to these two maps, the following lemma describes the special
fibre \(X_0\) and the geometric generic fibre \(X_\Omega\).

\begin{lemma}
\label{lem:generic-isocrystal}
Let \((Q,\phi)\) be a Breuil window with associated \(p\)-divisible
group \(X/R\).  Let \(L\) be a perfect field of characteristic \(p\),
and let \(\kappa\colon\Sfr\to W(L)\) be a ring homomorphism with
\(\kappa\circ\sigma=F\circ\kappa\).  Write \(\bar\kappa\colon R\to L\)
for its reduction modulo \(p\), and \(X_L\) for the base change of
\(X\) along \(\bar\kappa\).  Then there is an isomorphism of
\(F\)-isocrystals over \(W(L)[1/p]\)
\[
 \mathbb D(X_L)_{\Q}
 \simeq\bigl(Q\otimes_{\Sfr,\kappa}W(L)\bigr)[1/p],
 \qquad F^\sharp=p\phi^{-1},
\]
natural in the window \((Q,\phi)\).  In particular, it is compatible
with every endomorphism of the window.
\end{lemma}

\begin{proof}
We use the covariant Dieudonn\'e crystal \(\mathbb D\), as in
\cite{Lau2014}.  If \(S\) is a \(p\)-adically complete ring and
\(S\to R'\) is a surjection onto a ring with \(pR'=0\) whose kernel
carries divided powers compatible with those of \(p\), put
\(\mathbb D(G)_{S/R'}=\varprojlim_m\mathbb D(G)_{(S/p^mS\to R')}\).

\smallskip\noindent\emph{Step 0. Two frames.}
Put \(\mathcal W_L=(W(L),pW(L),L,F,p^{-1}F)\).  Since \(L\) is perfect,
\(pW(L)=VW(L)\) is the kernel of \(W(L)\to L\).  The frames
\(\mathcal B\) and \(\mathcal W_L\) satisfy the hypotheses at the
beginning of \cite[\S3C]{Lau2014}: all rings involved are
\(p\)-adically complete, \(\Sfr\) and \(W(L)\) are \(p\)-torsion free,
the ideals \(p\Sfr\) and \(pW(L)\) carry their canonical divided
powers, and \(\sigma=p\sigma_1\), respectively \(F=p\cdot p^{-1}F\), on
these ideals.  By \cite[Construction~3.16 and
Proposition~3.17]{Lau2014} we obtain crystalline windows
\(\mathcal P_{\mathcal B}(X)\) and \(\mathcal P_{\mathcal W_L}(X_L)\),
with underlying modules \(\mathbb D(X)_{\Sfr/R}\) and
\(\mathbb D(X_L)_{W(L)/L}\).  Their Frobenius is the crystalline
Frobenius computed through \(\sigma\), respectively \(F\), and their
divided Frobenius is determined by \(pF_1=F\) on the filtration,
because the base rings are \(p\)-torsion free.

\smallskip\noindent\emph{Step 1. The window over \(\Sfr\).}
Let $\mathcal P_0=(P_0,Q,F_0,F_1)$ be the $\mathcal B$-window
attached to $(Q,\phi)$; thus $P_0$ is its underlying
$\Sfr$-module.  We show that
\(\mathcal P_{\mathcal B}(X)\simeq\mathcal P_0\), naturally in \(X\).
The following argument is the \(E=p\) specialization of the proof of
\cite[Proposition~7.1]{Lau2014}.  We do not invoke that proposition
itself, since it is stated there under the additional assumption that
\(R\) has characteristic zero.  In the present case the
\(p\)-adically completed divided-power envelope of \(p\Sfr\) is
\(\Sfr\) itself: the canonical divided powers are
\(\gamma_n(px)=p^n x^n/n!\), and \(p^n/n!\in\Z_{(p)}\).  Moreover,
\(\sigma(p)/p=1\).  Thus the comparison homomorphism occurring in the
proof of loc.\ cit.\ is the identity, and the argument simplifies as
follows.

Write \(A=\mathbb W(R)=\varprojlim_n\mathbb W(R_n)\), where
\(R_n=R/(t^n)\).  Each \(R_n\) is admissible, and \(R=k[[t]]\) is an
admissible topological ring in the sense of
\cite[Definition~1.2]{Lau2014}.  The ring \(A\) is \(p\)-adically complete by
\cite[Proposition~1.14]{Lau2014} and \(p\)-torsion free, since
\(A\subset W(R)\subset W(\Omega)\).  Let \(\mathcal F_R\) denote the
Dieudonn\'e frame of \cite{Lau2014} with underlying ring \(A\).
Because \(pR=0\), the coefficient map is a strict frame homomorphism
\(\varkappa\colon\mathcal B\to\mathcal F_R\).  It is also compatible
with divided powers.  Let \(\delta\) denote the divided powers on the
augmentation ideal of \(A\) used in the Dieudonn\'e frame.  Since \(A\) is a \(p\)-torsion free \(\Z_p\)-algebra, every
 divided-power structure on an ideal of \(A\) satisfies
\(\delta_n(y)=y^n/n!\) in \(A[1/p]\).  Hence, for
\(x\in\Sfr\),
\[
 \delta_n\bigl(\varkappa(px)\bigr)
 =\frac{p^n\varkappa(x)^n}{n!}
 =\varkappa\!\left(\frac{p^n x^n}{n!}\right)
 =\varkappa\bigl(\gamma_n(px)\bigr).
\]
Both sides lie in $A$, so the equality in $A[1/p]$ holds already
in $A$.  The divided powers on $\mathbb W(R_n)$ are those of
\cite[\S1G, Lemma~1.17]{Lau2014}.  For $p=2$, this lemma applies
to the identity divided-power extension $R_n\to R_n$, since
$2R_n=0$.  The inverse-limit construction of $\mathbb W(R)$ in
\cite[\S1E]{Lau2014} shows that the projections
$A\to\mathbb W(R_n)$ respect divided powers.  Hence, for every $n$,
crystal base change along the divided-power morphism
\((\Sfr\to R)\to(\mathbb W(R_n)\to R_n)\) gives
\[
 \mathbb D(X)_{\Sfr/R}\otimes_{\Sfr}\mathbb W(R_n)
 \simeq\mathbb D(X_{R_n})_{\mathbb W(R_n)/R_n},
\]
compatibly with Hodge filtrations, Frobenius and the transition maps.
Here the filtration on a base-changed window is the inverse image of
the base-changed Hodge filtration modulo the augmentation ideal;
it agrees with the crystalline Hodge filtration because the Lie
module commutes with base change.
By \cite[Theorem~3.19]{Lau2014} for odd \(p\), and by
\cite[Corollary~3.24]{Lau2014} for \(p=2\), where \(2R_n=0\), the right
side, with its filtration and Frobenius, underlies the Dieudonn\'e
display \(\Phi_{R_n}(X_{R_n})\).  Passing to the inverse limit, as in
\cite[Lemma~2.16]{Lau2014}, identifies
\(\varkappa^*\mathcal P_{\mathcal B}(X)\) with \(\Phi_R(X)\) as
filtered \(F\)-modules, and then as windows, since \(A\) is
\(p\)-torsion free and \(F_1\) is therefore determined by \(pF_1=F\).
By the construction of the equivalence in
\cite[Corollary~6.7]{Lau2014}, \(\varkappa^*\mathcal P_0\simeq\Phi_R(X)\).
Since \(\varkappa\) is crystalline \cite[Theorem~6.6]{Lau2014}, the
functor \(\varkappa^*\) is fully faithful, so
\(\mathcal P_{\mathcal B}(X)\simeq\mathcal P_0\).  All isomorphisms
used are natural in \(X\).

\smallskip\noindent\emph{Step 2. Base change to \(W(L)\).}
Since \(\kappa(p\Sfr)\subseteq pW(L)\) and \(\kappa\sigma=F\kappa\),
\[
 p^{-1}F\bigl(\kappa(px)\bigr)=F(\kappa(x))=\kappa(\sigma(x))
 =\kappa(\sigma_1(px))\qquad(x\in\Sfr),
\]
so \(\kappa\colon\mathcal B\to\mathcal W_L\) is a strict frame
homomorphism.  It respects the canonical divided powers, since on
both \(p\Sfr\) and \(pW(L)\) they are
\(\gamma_n(px)=p^n x^n/n!\), and \(\kappa(p)=p\).  Hence, for each
\(m\), \(\kappa\)
induces a morphism of divided-power thickenings
\[
 \bigl(\Spec L\hookrightarrow\Spec W_m(L)\bigr)\longrightarrow
 \bigl(\Spec R\hookrightarrow\Spec\Sfr/p^m\Sfr\bigr)
\]
over \(\Spec\bar\kappa\).  These morphisms are compatible with the
transition maps in \(m\) and, since \(\kappa\sigma=F\kappa\), with the
Frobenius lifts.  By the base-change functoriality of the Dieudonn\'e crystal for
\(p\)-divisible groups \cite[\S3.3]{BBM}, the crystal property gives natural
isomorphisms
\[
 \mathbb D(X)_{(\Sfr/p^m\to R)}\otimes_{\Sfr/p^m}W_m(L)
 \simeq\mathbb D(X_L)_{(W_m(L)\to L)},
\]
compatible in \(m\) by the cocycle condition.  As
\(\mathbb D(X)_{\Sfr/R}\) is finite free, their inverse limit is an
isomorphism
\[
 \mathbb D(X)_{\Sfr/R}\otimes_{\Sfr,\kappa}W(L)\simeq
 \mathbb D(X_L)_{W(L)/L}.
\]
The Frobenius on each side is obtained from the Frobenius of the
crystal by evaluation along \(\sigma\), respectively \(F\).  Because
the relevant squares of divided-power morphisms commute, this
isomorphism is Frobenius-equivariant.  Thus
\(\kappa_*\mathcal P_{\mathcal B}(X)\simeq\mathcal P_{\mathcal W_L}(X_L)\)
as \(F\)-modules, naturally in \(X\).

\smallskip\noindent\emph{Step 3. Perfect fields.}
For the perfect field \(L\), the classical comparison of crystalline
and Dieudonn\'e theory \cite[\S4.2]{BBM} identifies
\(\mathbb D(X_L)_{W(L)/L}\) with the covariant Dieudonn\'e module
of \(X_L\); see also
\cite{BerthelotParfait}.  Its rational module is therefore the
isocrystal of \(X_L\): its slopes are the slopes of \(X_L\), and for
\(L\) algebraically closed its endomorphism algebra is \(\Endo(X_L)\)
\cite[Ch.~IV]{Demazure}.

\smallskip\noindent\emph{Step 4. Normalization.}
For the inclusion \(i\colon Q\hookrightarrow P_0\), the description
preceding \cite[Corollary~6.7]{Lau2014} gives $i=F_1^\sharp\phi$ and
$F_0^\sharp i^{(\sigma)}=pF_1^\sharp$.  Hence
$i^{-1}F_0^\sharp i^{(\sigma)}=p\phi^{-1}$ on $Q[1/p]$.  Combining
this with Steps~1--3 and inverting \(p\) gives the asserted
isomorphism.  This is the same normalization as in
\cite[Proposition~7.1]{Lau2014}: for \(E=p\), the Frobenius
\(x\mapsto1\otimes\phi^{-1}(Ex)\) on \(Q^{(\sigma)}\) is carried to
\(p\phi^{-1}\) on \(Q\) by the Frobenius-equivariant isogeny
\(p\phi^{-1}\colon Q^{(\sigma)}\to Q\).

Finally, every isomorphism above is natural in \(X\), and endomorphisms
of the window are endomorphisms of \(X\) by
\cite[Corollary~6.7]{Lau2014}.  This gives compatibility with
endomorphisms of the window.
\end{proof}

\begin{lemma}\label{lem:val}
Put \(K_0=k((t))\) and let \(\mathcal C_{K_0}\) be the \(p\)-adic
completion of \(\Sfr[1/t]\).  It is a complete discrete valuation
ring with uniformizer \(p\) and residue field \(K_0\).  For a nonzero
\(x=\sum_i a_it^i\in\Sfr\) one has
\[
 v_p(x)=\min_i v_p(a_i).
\]
In particular, \(t\), \(1-t\), and every series with a unit
coefficient are units in \(\mathcal C_{K_0}\).
\end{lemma}

\begin{proof}
The ring \(\mathcal C_{K_0}\) is \(p\)-adically complete and separated,
\(p\)-torsion free, and \(\mathcal C_{K_0}/p=K_0\).  Every element
with nonzero reduction is therefore a unit.  If \(0\ne z\in
\mathcal C_{K_0}\), separatedness gives a largest integer \(m\ge0\)
such that \(z\in p^m\mathcal C_{K_0}\).  Thus \(z=p^m u\), where the
reduction of \(u\) modulo \(p\) is nonzero, so \(u\) is a unit.  Hence
every nonzero element is a power of \(p\) times a unit, and
\(\mathcal C_{K_0}\) is a discrete valuation ring with uniformizer
\(p\).  After dividing \(x\) by \(p^{\min_i v_p(a_i)}\), its
reduction is a nonzero element of \(k[[t]]\subset K_0\).
\end{proof}

\begin{lemma}
\label{lem:generic-valuations}
Extend \(\sigma\) to the ring \(\mathcal C_{K_0}\) of
\cref{lem:val} by \(\sigma(t^{-1})=t^{-p}\).  The map \(\kappa_\eta\)
extends uniquely to a continuous ring homomorphism
\(\mathcal C_{K_0}\to W(\Omega)\).  The extension satisfies
\(\kappa_\eta\sigma=F\kappa_\eta\) and preserves \(p\)-adic valuations.
Let \(T\in GL_n(\Sfr[1/p])\), put
\(T_N=T\,\sigma(T)\cdots\sigma^{N-1}(T)\), and let
\(H=\kappa_\eta(T)\,F\) act on \(W(\Omega)[1/p]^n\), with
\(\Lambda=W(\Omega)^n\).  Then \(H^N(\Lambda)=\kappa_\eta(T_N)\Lambda\).
Using the convention \(v_p(0)=+\infty\), we consequently have
\[
 \max\{c\in\Z:H^N(\Lambda)\subseteq p^c\Lambda\}
 =\min_{i,j}v_p\bigl((T_N)_{ij}\bigr),
 \qquad
 v_p\bigl(\det\kappa_\eta(T)\bigr)=v_p(\det T),
\]
where \(v_p\) on the right is the valuation of \(\mathcal C_{K_0}\)
computed in \cref{lem:val}.
\end{lemma}

\begin{proof}
Since \([t]\) is a unit of \(W(\Omega)\), \(\kappa_\eta\) extends to
\(\Sfr[1/t]\), and then to its \(p\)-adic completion, because
\(W(\Omega)\) is \(p\)-adically complete.  Uniqueness and Frobenius
compatibility follow by continuity.  By \cref{lem:val}, the extension
is a local homomorphism of discrete valuation rings with the
common uniformizer \(p\), inducing \(K_0\hookrightarrow\Omega\) on
residue fields; hence it preserves valuations.

By \(\kappa_\eta\sigma=F\kappa_\eta\), for \(x\in W(\Omega)[1/p]^n\),
\[
 H^N(x)=\kappa_\eta(T)\,F(\kappa_\eta(T))\cdots
 F^{N-1}(\kappa_\eta(T))\,F^N(x)
 =\kappa_\eta(T_N)\,F^N(x).
\]
Since \(\Omega\) is perfect, \(F^N\) is bijective on \(\Lambda\), so
\(H^N(\Lambda)=\kappa_\eta(T_N)\Lambda\).  A matrix \(B\) satisfies
\(B\Lambda\subseteq p^c\Lambda\) exactly when all its entries have
valuation at least \(c\).  Valuation preservation gives the first
formula; the second holds because \(\kappa_\eta\) is a ring
homomorphism.
\end{proof}

\begin{lemma}\label{lem:growth}
Let \((N,H)\) be a nonzero isocrystal over a perfect field and let
\(\Lambda\subset N\) be a lattice.  Then
\[
 \lambda_{\min}(H)=\lim_{n\to\infty}\frac1n
 \max\{c\in\Z:H^n(\Lambda)\subseteq p^c\Lambda\}.
\]
\end{lemma}

\begin{proof}
Both sides are unchanged by extension to an algebraic closure, using
faithful flatness for the lattice condition.
On a standard Dieudonn\'e--Manin lattice of slope \(c/e\), one has
\(H^e\Lambda=p^c\Lambda\); thus the displayed maximum is
\(nc/e+O(1)\).  For a direct sum of standard lattices it is the minimum
of these quantities.  Any other lattice is commensurable with this
one, which changes the maximum by a bounded amount independent of
\(n\).  See \cite[Ch.~IV]{Demazure}.
\end{proof}

\begin{lemma}\label{lem:smith}
Let \(A\in M_N(\Sfr)\) and let \(\bar A\) denote its reduction
modulo \(p\).  Suppose \(\bar A\) has rank \(c\) over \(k((t))\),
some \(c\times c\) minor of \(\bar A\) is a unit of \(R\), and
\(\det A=p^{N-c}\varepsilon\) with \(\varepsilon\in\Sfr^\times\).  Then
\[
 A=U\operatorname{diag}(I_c,pI_{N-c})V
 \quad\text{for some }U,V\in GL_N(\Sfr).
\]
Consequently \(pA^{-1}\) is integral,
\(\coker(A)\simeq R^{N-c}\), and
\(\coker(pA^{-1})\simeq R^c\).
\end{lemma}

\begin{proof}
Move the unit minor to the upper-left block.  Its lift \(A_{11}\)
has determinant outside the maximal ideal \((p,t)\) of \(\Sfr\), and is
therefore invertible over \(\Sfr\).  Block operations reduce
\(A\) to \(\operatorname{diag}(A_{11},A^\circ)\), with
\(A^\circ=A_{22}-A_{21}A_{11}^{-1}A_{12}\).  The rank hypothesis gives
\(\bar A^\circ=0\) over \(k((t))\), hence over \(R\), so
\(A^\circ=pA'\).  The determinant hypothesis makes \(\det A'\) a unit.
This proves the normal form and both cokernel assertions, including
the cases of empty blocks.
\end{proof}

\subsection{Cyclic actions and reduced centralizers}\label{sec:reduced-centralizers}

Let $\Delta_h/\Q_p$ be the central division algebra of invariant $1/h$.
We use its standard maximal order
\[
 \mathcal O_h=W(\F_{p^h})[\Pi],\qquad
 \Pi^h=p,\qquad \Pi\alpha=\sigma(\alpha)\Pi.
\]
For a rational Dieudonn\'e module $N$ over a perfect field containing
$\F_{p^h}$, with $\Delta_h$-action, the
unramified subfield gives a grading $N=\bigoplus_{i\in\Z/h}N_i$,
indexed so that $F:N_i\to N_{i+1}$ and $\Pi:N_i\to N_{i-1}$.
The \emph{reduced operator} $G=\Pi F$ preserves each $N_i$ and
satisfies $G^h=pF^h$.  Therefore the relation between a reduced slope
$\nu$ on $N_0$ and its full covariant slope $\lambda$ is
\begin{equation}\label{eq:boxed-slope-convention}
 \nu=\lambda+\frac1h,\qquad
 \nu^{[m]}\text{ on }N_0
 \quad\Longleftrightarrow\quad
 \lambda^{[hm]}\text{ on }N.
\end{equation}
The multiplicity factor comes from the $h$ graded pieces, which are
identified by $\Pi$ after inverting $p$.

\begin{proposition}\label{prop:rank2}
The restriction map from $\End_{\Delta_h}(N)$ to $\End(N_0,G)$ is an
isomorphism.  If $X$ is an $\mathcal O_h$-linear $p$-divisible group,
then the
standard matrix action on $X^m$ gives
\[
 \Endo_{M_m(\Delta_h)}(X^m)\simeq\Endo_{\Delta_h}(X).
\]
In particular, if the geometric reduced isocrystal of $X$ has exactly
two slopes, with denominators $r$ and $s$ in lowest terms and
multiplicities $r$ and $s$, then the corresponding geometric
centralizers are
$\mathcal D_r\times\mathcal D_s$, with central division factors over
$\Q_p$ of dimensions $r^2,s^2$.
\end{proposition}

\begin{proof}
A $\Delta_h$-linear endomorphism commutes with $\Pi$ and $F$, hence
with $G$, and is determined by its restriction to $N_0$.  Conversely,
an endomorphism $c_0$ commuting with $G$ extends to $N_i$ as
$c_i=\Pi^{-i}c_0\Pi^i$.  The relation $\Pi^h=p$ makes this independent
of the representative of $i$.  Since $G$ commutes with $\Pi$, the maps
$c_i$ commute with $G$ on every graded piece.  Moreover, on $N_i$ the
unramified subfield acts through the embedding $\sigma^i$, so the $c_i$
commute with this action.  Hence the resulting endomorphism commutes
with the unramified subfield, with $\Pi$, and with $F=\Pi^{-1}G$.

For the matrix action, commuting with the diagonal matrix units makes
an endomorphism of $X^m$ diagonal, and commuting with the off-diagonal
units makes its diagonal entries equal; the remaining condition is
$\Delta_h$-linearity.  Finally, Dieudonn\'e--Manin theory makes each
reduced slope piece simple when its rank equals its denominator, gives
no morphisms between the two slopes, and identifies their endomorphism
algebras as central division algebras of the stated dimensions
\cite[Ch.~IV]{Demazure}.
\end{proof}

\subsection{The matrix and two lattice chains}

Let $k$ be algebraically closed of characteristic $p$, put
$R=k[[t]]$ and $\Sfr=W(k)[[t]]$, and let $s>r\ge2$ and $n=r+s$.
The matrix $M_{s,r}(t)$ defined in \eqref{eq:gen-edges} below
interpolates linearly between two monomial matrices:
\[
 M_{s,r}(t)=(1-t)\,M_{s,r}(0)+t\,M_{s,r}(1).
\]
Their cycle types are $(r-1,s+1)$ and $(r,s)$: the first gives the
isoclinic special fibre, while the second has the desired generic
Newton polygon.
The operator has two integral realizations: one with an action of a
division algebra of degree $s$, used for the slope bound and the
quadratic-CM construction, and one with a quaternion action, used for
the main theorem.
Set $u=r-1$, $x=r$, $v=r+1$, and $w=r+2$.
An edge $j\to z$ labelled by $\epsilon$ means that the $(z,j)$-entry
of $M=M_{s,r}(t)$ is $\epsilon$.  The nonzero entries of $M$ are
\begin{equation}\label{eq:gen-edges}
\begin{gathered}
 j\longrightarrow j+1:p
 \quad(1\le j\le r-2\text{ or }w\le j\le n-1),
 \qquad n\longrightarrow x:p^2,\\
 \begin{array}{c|ccc}
  &u&x&v\\ \hline
  1&p(1-t)&0&pt\\
  v&t&1-t&0\\
  w&0&pt&p(1-t)
 \end{array}
 \quad\text{(rows are targets, columns are sources).}
\end{gathered}
\end{equation}
Write
\begin{equation}\label{eq:gen-q}
 \vartheta(t)=(1-t)^3+t^3=3t^2-3t+1\in\Sfr^\times.
\end{equation}

\begin{lemma}\label{lem:gen-det}
One has $\det M_{s,r}(t)=(-1)^n p^n \vartheta(t)$.
\end{lemma}

\begin{proof}
In a cycle cover all edges outside the sources $u,x,v$ are forced,
with total $p$-weight $n-2$.  The remaining sources map either to
$(1,v,w)$ or to $(v,w,1)$, contributing respectively
$p^2(1-t)^3$ and $p^2t^3$.  The resulting covers have cycle lengths
$(r-1,s+1)$ and $(r,s)$, so both signs are $(-1)^{n-2}=(-1)^n$.
For $r=2$ the first cycle is a loop.
\end{proof}

Define the levels and their two-level compression by
\begin{equation}\label{eq:gen-level}
 L(j)=
 \begin{cases}
 j,&1\le j\le r-1,\\
 1,&j=r,\\
 r,&j=r+1,\\
 j-r,&j\ge r+2,
 \end{cases}
 \qquad
 \ell_s(j)=L(j),\qquad
 \ell_2(j)=1+\mathbf1_{L(j)\ge r}.
\end{equation}
For either $h=s$ or $h=2$, use $\ell=\ell_h$ and set
\[
 H_i=\operatorname{diag}\bigl(p^{\mathbf1_{\ell(j)\le i}}\bigr)_{j=1}^n
 \quad(0\le i\le h),\qquad H_{i+h}=pH_i,
 \qquad A_i=H_{i+1}^{-1}MH_i.
\]
Thus $H_0=I_n$, $H_h=pI_n$, and $A_i$ is periodic of period $h$.
The successive colengths $d_i=\#\{j:\ell(j)=i+1\}$ are
\begin{equation}\label{eq:gen-signature}
 \begin{array}{ll}
 h=s:&d_i=2\ (0\le i<r),\quad d_i=1\ (r\le i<s),\\
 h=2:&d_0=2r-2,\quad d_1=s-r+2.
 \end{array}
 \qquad \sum_i d_i=n.
\end{equation}

\begin{proposition}\label{prop:gen-elem}
For both lattice chains and $0\le i<h$, $A_i$ is integral, its reduction modulo $p$
has rank $d_i$ over $k((t))$ and a unit $d_i\times d_i$ minor over $R$, and
$\det A_i=(-1)^n p^{n-d_i}\vartheta(t)$.  Consequently
\[
 A_i\sim\operatorname{diag}(I_{d_i},pI_{n-d_i}),\qquad
 \Phi_i=pA_i^{-1}\in M_n(\Sfr),\qquad
 \coker\Phi_i\simeq R^{d_i}.
\]
\end{proposition}

\begin{proof}
For an edge $j\to z$ of $p$-weight $m$, its entry in $A_i$ has exponent
$m+\mathbf1_{\ell(j)\le i}-\mathbf1_{\ell(z)\le i+1}$.
Only the weight-zero edges $u\to v$ and $x\to v$ could give a negative
exponent.  For $h=s$ their levels are $(r-1,r)$ and $(1,r)$;
for $h=2$ both go from level $1$ to level $2$.
The exponent is nonnegative in every case, proving integrality.

For $h=s$, reduction modulo $p$ keeps exactly the entries of exponent
zero: $u\to v$ for every $i$; $x\to v$ for $i=0$ or $i\ge r-1$;
$u\to1$ for $0\le i\le r-2$; $v\to1$ for $0\le i\le r-1$; and $v\to w$
for $1\le i\le r-1$.  Every other edge, including both chains,
$x\to w$ and $n\to x$, has exponent one.  The following minors of
$A_i\bmod p$ have determinant $(1-t)^2$:
\[
 \begin{array}{c|c|c}
 i&\text{rows}&\text{columns}\\ \hline
 0&(1,v)&(u,x)\\
 1\le i\le r-2&(1,w)&(u,v)\\
 r-1&(v,w)&(x,v).
 \end{array}
\]
The support has rank at most two: for $i=0$ only rows $1,v$ survive;
for $1\le i\le r-2$ only columns $u,v$ survive; and for $i=r-1$
the nonzero rows are $1,v,w$, with rows $1,w$ supported only in column
$v$.  For $i\ge r$ only row $v$ survives and contains the unit $1-t$.
This proves the claimed ranks and unit minors.  The middle range is
empty when $r=2$.

For $h=2$, the first level is
$\{1,\ldots,r\}\cup\{r+2,\ldots,2r-1\}$, of size $d_0=2r-2$;
the complementary level has size $d_1=s-r+2$.
For $i=0$ an edge $j\to z$ of weight $m$ has exponent
$m-\mathbf1_{\ell(z)=1}$.  The surviving entries of $A_0\bmod p$ are
therefore $u\to v$, $x\to v$, and the weight-one edges whose target
lies in the first level; the edge $n\to x$ keeps exponent one.  Their
targets, hence the nonzero rows, are exactly
$\{1,\ldots,r-1\}\cup\{r+1,\ldots,2r-1\}$.
At $t=0$, unit entries form the cycle on $1,\ldots,r-1$ and the chain
$r\to r+1\to\cdots\to2r-1$.  Since $M(0)$ is monomial, these
$d_0$ distinct sources and targets give a $d_0\times d_0$ minor
with value $\pm1$ at $t=0$.  Hence its determinant is a unit in $R$.
For $r=2$ these are the loop at $1$ and the edge $2\to3$;
in particular, the first level is exactly $\{1,2\}$ and the same
rank formula holds.
For $i=1$ the exponent is $m+\mathbf1_{\ell(j)=1}-1$.  The surviving
entries of $A_1\bmod p$ are therefore $u\to v$ and $x\to v$, the edges
$v\to1$ and $v\to w$ leaving the second-level vertex $v$, and the
chain $2r\to2r+1\to\cdots\to n$; again $n\to x$ keeps exponent one.
The block on rows $(1,v,w)$ and columns $(u,x,v)$ is
\[
 \begin{pmatrix}
 0&0&t\\ t&1-t&0\\ 0&0&1-t
 \end{pmatrix}.
\]
It has rank two and the unit minor $(1-t)^2$ on its last two rows and
columns.  The chain $2r\to2r+1\to\cdots\to n$ has unit entries and
rows and columns disjoint from this block.  It contributes $n-2r$ to the rank and extends
that minor to a unit minor of size $2+n-2r=d_1$.

Finally, $\det H_{i+1}=p^{d_i}\det H_i$.  Hence the definition
$A_i=H_{i+1}^{-1}MH_i$ and \cref{lem:gen-det} give
\[
 \det A_i=(-1)^n p^{n-d_i}\vartheta(t).
\]
Since $\vartheta(0)=1$, \cref{lem:smith} applies to both chains and
gives the elementary divisors and the assertions about $\Phi_i$.
\end{proof}

\subsection{The associated windows}

\begin{theorem}\label{thm:gen-window}
For $h=s$ and $h=2$, put $Q_i=\Sfr^n$ and
$\cQ=\bigoplus_{i\in\Z/h\Z}Q_i$.  The maps
$\phi|_{Q_{i+1}}=\Phi_i:Q_{i+1}\to Q_i^{(\sigma)}$ define an
$\cO_h$-linear Breuil window of height $nh$ and dimension $n$.
Its associated $p$-divisible group $X_h/R$ has isoclinic special fibre
of covariant slope $1-1/h$.  On the graded piece $Q_0$ the reduced
Frobenius is $G=M_{s,r}(t)\sigma$, independently of $h$.
\end{theorem}

\begin{proof}
By \cref{prop:gen-elem}, $\phi$ is injective with cokernel annihilated
by $p$ and free over $R$ of rank $\sum_i d_i=n$.
For $\alpha\in W(\F_{p^h})$, let $\varphi_\alpha$ act on $Q_i$
as $\sigma^i(\alpha)$.
Its scalar on $Q_i^{(\sigma)}$ is $\sigma^{i+1}(\alpha)$, the same
as on $Q_{i+1}$.
Set $C_i=H_iH_{i-1}^{-1}$ and let $\Pi:Q_i\to Q_{i-1}$ have matrix
$C_i$.  The diagonal matrices $C_i$ have entries $1,p$, are fixed by
$\sigma$, and satisfy
\[
 C_h\cdots C_1=pI_n,\qquad
 C_{i+1}A_i=A_{i-1}C_i,\qquad
 C_i\Phi_i=\Phi_{i-1}C_{i+1}.
\]
Since the $C_i$ are diagonal and periodic, every product of $h$
consecutive $C_i$ equals $pI_n$.  Thus $\Pi$ is an integral window
endomorphism with $\Pi^h=p$.
The chosen weights give
$\Pi\varphi_\alpha=\varphi_{\sigma(\alpha)}\Pi$, so this is an
$\cO_h$-action.  Since $\Phi_i=pA_i^{-1}$, the covariant rational
Frobenius $F=p\phi^{-1}$ restricts to $A_i\sigma$ on $Q_i[1/p]$.
It commutes with $\Pi$; on $Q_0[1/p]$ one has
$G=\Pi F=C_1A_0\sigma=M\sigma$.
At $t=0$, $M$ is monomial with cycles of lengths $r-1,s+1$ and
weights $r-1,s+1$, so its reduced slopes all equal $1$.
Since $G^h=pF^h$, a reduced slope $\nu$ of multiplicity $m$ gives
full slope $\nu-1/h$ of multiplicity $hm$.
The equivalence of \cite[Corollary~6.7]{Lau2014} gives $X_h$.
\Cref{lem:generic-isocrystal}, applied to $\kappa_0$ and to
$\kappa_\eta$, identifies its special and geometric generic
isocrystals, $\cO_h$-equivariantly, with those defined by these
matrices.  In particular, the special fibre is isoclinic of
covariant slope $1-1/h$.
\end{proof}

\subsection{The generic Newton polygon}

After the scalar extension of \cref{lem:generic-valuations}, $t$, $1-t$,
and all their Frobenius transforms are $p$-adic units: their reductions
are $t^{p^j}$ and $1-t^{p^j}$ in $k((t))$.
Thus the valuation of each term in an iterated Frobenius matrix is the
weight of its directed walk.  The directed graph of the nonzero entries
of $M$ has exactly four simple cycles:
\begin{equation}\label{eq:gen-cycles}
 \begin{array}{c|l|c|c|c}
 &\text{cycle}&\text{length}&\text{weight}&\text{mean}\\ \hline
 L_0&(1,\ldots,u,1)&r-1&r-1&1\\
 L_1&(1,\ldots,u,v,1)&r&r-1&1-1/r\\
 R_0&(v,w,\ldots,n,x,v)&s+1&s+1&1\\
 R_1&(w,\ldots,n,x,w)&s&s+1&1+1/s.
 \end{array}
\end{equation}
Indeed, a simple cycle meeting $1$ follows the forced chain to $u$ and
closes either by $u\to1$ or by $u\to v\to1$.  Taking $v\to w$
instead would force a repetition of $v$ or $w$ before returning to $1$.
A cycle avoiding $1$ cannot meet $u$; according as it meets $v$ or not,
the forced remaining edges give $R_0$ or $R_1$.

\begin{lemma}\label{lem:gen-min}
The smallest slope of the geometric generic reduced operator
$M_{s,r}(t)\sigma$ is $1-1/r$.
\end{lemma}

\begin{proof}
Put $\alpha=(r-1)/r$.  Delete simple closed subwalks from a walk
$\gamma$ of length $N$ until a simple path remains.  This path has
length at most $n-1$ and nonnegative weight; every removed cycle has
mean weight at least $\alpha$.  Hence
$w(\gamma)\ge\alpha(N-n+1)$.
If $c_N$ is the minimum entry valuation of
$M\sigma(M)\cdots\sigma^{N-1}(M)$, then
$\liminf_N c_N/N\ge\alpha$.

Fix a vertex of $L_1$.  The walk $L_1^m$ has length $mr$ and weight
$m(r-1)$, and is the unique walk with these endpoints, length, and
weight.  Indeed, any closed walk of mean $\alpha$ decomposes entirely
into simple cycles of mean $\alpha$; by \eqref{eq:gen-cycles} each must
be $L_1$, whose orientation and initial vertex determine the walk.
Consequently the corresponding diagonal entry of
$M\sigma(M)\cdots\sigma^{mr-1}(M)$ has a unique term of valuation
$m(r-1)$, all other terms having larger valuation.  There is no
cancellation at this valuation, so $c_{mr}\le m(r-1)$.
By \cref{lem:generic-valuations}, these valuations are preserved by the
scalar extension, and \cref{lem:growth} identifies
$\lim_N c_N/N$ with the smallest slope.  The lower bound and
this subsequence give the assertion.
\end{proof}

\begin{theorem}\label{thm:gen-polygon}
For $s>r\ge2$, the geometric generic reduced operator $M_{s,r}(t)\sigma$
has slopes
\[
 \left(1-\frac1r\right)^{[r]}
 \cup\left(1+\frac1s\right)^{[s]}.
\]
The full covariant slopes of the geometric generic fibre of
the degree-$s$ group $X_s$ are
\[
 \left(1-\frac1r-\frac1s\right)^{[rs]}\cup1^{[s^2]}.
\]
\end{theorem}

\begin{proof}
The determinant formula gives total reduced slope $n$, hence average
$1$.  For the degree-$s$ lattice chain, both $A_i$ and $pA_i^{-1}$
are integral.  After passage to the perfect geometric generic field,
the lattice is stable under $F$ and $V=pF^{-1}$, so every full slope
satisfies $0\le\mu\le1$.  Since $G^s=pF^s$, every reduced slope
satisfies $\nu\le1+1/s$.

By \cref{lem:gen-min}, the smallest slope is $(r-1)/r$.
Its denominator divides its multiplicity; write that multiplicity as
$mr$, with $m\ge1$.  This slope cannot have multiplicity $n$, since it
is less than the average $1$.  The remaining slopes have average
\[
 1+\frac{m}{n-mr}\le1+\frac1s,
 \qquad\text{hence}\qquad m(r+s)\le n=r+s.
\]
Thus $m=1$.  The remaining $s$ slopes have average $1+1/s$ and are
bounded above by that number, so all equal it.  Subtracting $1/s$ and
multiplying multiplicities by $s$ gives the full slopes of $X_s$.
\end{proof}

\begin{remark}
\label{rem:interpolation}
At $t=0$, the surviving edges form the cycles $L_0$ and $R_0$,
both of mean weight $1$.  At $t=1$, they form $L_1$ and $R_1$,
with the reduced slopes of \cref{thm:gen-polygon}.  Thus the family
changes the cycle type $(r-1,s+1)$ into $(r,s)$.  The two terms
$(1-t)^3$ and $t^3$ in $\vartheta(t)$ are precisely the contributions
of these two endpoint cycle covers.
\end{remark}

\begin{proposition}\label{prop:two-level-window-complete}
For every $s>r\ge2$, the degree-$2$ construction gives an
$\cO_2$-linear group $X_2/R$ of height $2(r+s)$ and dimension
$r+s$.  Its special fibre has covariant
slope $1/2$, and its full geometric generic slopes are
\[
 \left(\frac12-\frac1r\right)^{[2r]}
 \cup\left(\frac12+\frac1s\right)^{[2s]}.
\]
\end{proposition}

\begin{proof}
Combine the reduced slopes of \cref{thm:gen-polygon} with
\cref{thm:gen-window} for $h=2$, using the relation
$\lambda=\nu-1/2$.
\end{proof}

For $r=2$ the compression still applies, but the generic slope $0$
yields a slope $1$ in the dual factor of the polarized construction,
that is, positive $p$-rank (see \cref{cor:r2-variant}).  We therefore
take $r\ge3$ when constructing varieties of $p$-rank zero.
The degree-$s$ realization has a slope-$1$ component and will be used
for the smaller quadratic-CM construction with positive $p$-rank.

\section{Honda--Tate seeds}\label{sec:seeds}

To globalize the local deformations of Section~\ref{sec:explicit-matrix},
we need abelian varieties with compatible local endomorphism algebras.
Fix a prime $p$.  We use Honda--Tate theory in the following form.
The existence and classification statement is recalled in
\cite[Part~II, Theorem~2]{WaterhouseMilne}; for the endomorphism algebra
and the local invariant formula we use
\cite[Theorem~2 and (7), pp.~140--142]{Tate66}.  Let $\pi$ be a $q$-Weil
number, let $B/\F_q$ be a simple abelian variety with Weil number $\pi$,
and put $L=\Q(\pi)$ and $D=\Endo_{\F_q}(B)$.  Then
\begin{equation}\label{eq:HT-invariant-convention}
 \inv_v(D)=[L_v:\Q_p]\frac{v(\pi)}{v(q)}\quad(v\mid p),
 \qquad 2\dim B=[L:\Q]\deg_LD.
\end{equation}
The finite invariants away from $p$ vanish; real places have invariant
$1/2$.  The degree of $D$ is the least common denominator of its local
invariants: this follows from the local-global injection for the Brauer
group and equality of period and index over number fields
\cite[Theorem~1.5.36(i),(iv)]{Poonen}.
In the covariant convention used here, the slope at
$v\mid p$ is $1-v(\pi)/v(q)$.

We will repeatedly use the following geometric base-change criterion.
If $\Q(\pi^m)=L$ for every $m\ge1$, then $B$ is geometrically simple and
$\Endo(B_{\overline{\F}_p})=D$.  Indeed, the Frobenius of
$B_{\F_{q^m}}$ is $\pi^m$.  Since $[\Q(\pi^m):\Q]=[L:\Q]$, the $m$-th
powers of the conjugates of $\pi$ are distinct and are exactly the
conjugates of $\pi^m$.  Hence the characteristic polynomial of $\pi^m$
on $B$ is a power of the minimal polynomial of $\pi^m$, and
$B_{\F_{q^m}}$ is isogenous to a power $B_m^{e}$ of the simple variety
$B_m$ attached to $\pi^m$.  The centre and all local invariants
attached to $\pi^m$ coincide with those of $\pi$, because the ratios
$v(\pi^m)/v(q^m)$ are unchanged; so $\dim B_m=\dim B$ and $e=1$.
Thus $B_{\F_{q^m}}$ is simple, and the inclusion of its original
endomorphism algebra is an equality by dimension.  Abelian subvarieties
and endomorphisms over $\overline{\F}_p$ descend to finite extensions.
This proves the assertion; see
\cite[Part~I, Theorems~8--9]{WaterhouseMilne}.

\subsection{Quadratic seeds}

\begin{proposition}\label{prop:seeds}
For every integer $b\ge3$, put
\[
 \alpha=\frac{1+\sqrt{1-4p}}2,\qquad K_p=\Q(\alpha),\qquad
 \pi_{p,b}=\alpha\bar\alpha^{\,b-1},\qquad q=p^b.
\]
Then $p\mathcal O_{K_p}=\pfr\bar\pfr$, where $\pfr=(\alpha)$, and
$\pi_{p,b}$ is a $q$-Weil number satisfying
$\Q(\pi_{p,b}^m)=K_p$ for every $m\ge1$.  Its Honda--Tate variety
$B_{p,b}/\F_q$ is geometrically simple of dimension $b$ and $p$-rank
zero.  Its geometric endomorphism algebra $D_{p,b}$ is a central
division algebra of degree $b$ over $K_p$, with
\[
 \inv_{\pfr}(D_{p,b})=\frac1b,\qquad
 \inv_{\bar\pfr}(D_{p,b})=1-\frac1b,
 \qquad [D_{p,b}:\Q]=2b^2.
\]
All these endomorphisms are defined over $\F_q$.
\end{proposition}

\begin{proof}
The identities $\alpha+\bar\alpha=1$ and $\alpha\bar\alpha=p$ show
that $(\alpha)$ and $(\bar\alpha)$ are distinct prime ideals of norm
$p$.  Both complex conjugates of $\pi_{p,b}$ have absolute value
$p^{b/2}$, and
$(\pi_{p,b}^m)=\pfr^m\bar\pfr^{m(b-1)}$ is not fixed by conjugation.
Thus no positive power of $\pi_{p,b}$ is rational.  Since
$K_p/\Q$ is quadratic, it follows that $\Q(\pi_{p,b}^m)=K_p$ for every
$m\ge1$.  The invariant formula gives the displayed values, whose
denominator is $b$, and
the dimension formula gives $\dim B_{p,b}=b$.  The preceding
base-change criterion proves the geometric assertions.  The two
covariant slopes are $1-1/b$ at $\pfr$ and $1/b$ at $\bar\pfr$, each
with multiplicity $b$, so the $p$-rank is zero.  The argument applies
also to $p=2$ and does not require $1-4p$ to be squarefree.
\end{proof}

\subsection{A quartic CM field split at \texorpdfstring{$p$}{p}}

\begin{lemma}\label{lem:quartic-cm-split-complete}
Let $K_p=\Q(\sqrt{1-4p})$ be the field of \cref{prop:seeds}.  For odd
$p$ write $p^2+4=m_0c_0^2$ with $m_0$ squarefree and put
$K^+=\Q(\sqrt{m_0})$; for $p=2$ put $K^+=\Q(\sqrt{17})$.  Then
$K=K_pK^+$ is a quartic CM field in which $p$ splits completely.
\end{lemma}

\begin{proof}
For odd $p$, the integer $p^2+4$ lies strictly between $p^2$ and
$(p+1)^2$, so it is not a square and $m_0>1$.  Since
$p\nmid p^2+4=m_0c_0^2$, the prime $p$ divides neither $m_0$ nor $c_0$.
Reducing the equality $p^2+4=m_0c_0^2$ modulo $p$ gives
$m_0\equiv(2c_0^{-1})^2\pmod p$, a nonzero square.  Hence $p$
splits in the real quadratic field $K^+$; for $p=2$ this holds because
$17\equiv1\pmod8$.  By \cref{prop:seeds}, $p$ splits in the imaginary
quadratic field $K_p$.

The field $K=K^+(\sqrt{1-4p})$ is a totally imaginary quadratic
extension of the totally real field $K^+$, hence a quartic CM field,
biquadratic over $\Q$.  Restriction maps the decomposition group in
$\Gal(K/\Q)$ of a prime above $p$ onto the corresponding
decomposition groups in $\Gal(K_p/\Q)$ and $\Gal(K^+/\Q)$, which are
trivial.  Since $\Gal(K/\Q)\to\Gal(K_p/\Q)\times\Gal(K^+/\Q)$ is
injective, the decomposition group is trivial, so $p$ splits completely
in $K$.
\end{proof}

We use complete splitting through the following lemma.  For a prime
$\mathfrak P\mid p$ of a number field in which $p$ splits completely, we
normalize $v_{\mathfrak P}(p)=1$.

\begin{lemma}\label{lem:pair-valuations}
Let $K$ be a number field in which $p$ splits completely, and let
$x\in K^\times$.  If $\Q(x)\ne K$, then each value of the function
$\mathfrak P\mapsto v_{\mathfrak P}(x)$ on the primes $\mathfrak P\mid p$
of $K$ is attained at a number of primes divisible by
$[K:\Q(x)]\ge2$.  In particular, if some value is attained at exactly
one prime above $p$, then $\Q(x)=K$.
\end{lemma}

\begin{proof}
Put $F=\Q(x)$ and $t=[K:F]\ge2$.  The prime $p$ splits completely in
$F$, because ramification indices and residue degrees are multiplicative
in towers.  Each prime $\mathfrak q$ of $F$ above $p$ therefore lies
below exactly $t$ primes of $K$, all unramified over $\mathfrak q$, and
$v_{\mathfrak P}(x)=v_{\mathfrak q}(x)$ for each of them.  Thus the
primes of $K$ above $p$ are partitioned into sets of size $t$ on which
the valuation of $x$ is constant.
\end{proof}

Fix the field $K$ of \cref{lem:quartic-cm-split-complete}, and write
$p\mathcal O_K=\mathfrak p_1\bar{\mathfrak p}_1
\mathfrak p_2\bar{\mathfrak p}_2$, with the bar denoting CM conjugation.
This labelling is possible because each prime of $K^+$ above $p$ splits
in the quadratic extension $K/K^+$, whose nontrivial automorphism is CM
conjugation and permutes the two primes above it transitively.

\subsection{A quartic seed with different local indices}

\begin{proposition}\label{prop:quartic-seed-complete}
For every even integer $d\ge4$, there exist a finite field $\F_q$ of
characteristic $p$ and a geometrically simple abelian variety
$B_d/\F_q$ of dimension $2d$ and $p$-rank zero whose geometric
endomorphism algebra
$D_d=\Endo_{\F_q}(B_d)$ is a central division algebra of degree $d$ over $K$.
Its invariants at
$(\mathfrak p_1,\bar{\mathfrak p}_1,
  \mathfrak p_2,\bar{\mathfrak p}_2)$ are
\[
             \left(\frac12,\frac12,\frac1d,1-\frac1d\right),
             \qquad [D_d:\Q]=4d^2.
\]
\end{proposition}

\begin{proof}
Set
$\mathfrak a=\mathfrak p_1^{d/2}\bar{\mathfrak p}_1^{d/2}
\mathfrak p_2\bar{\mathfrak p}_2^{d-1}$, so
$\mathfrak a\bar{\mathfrak a}=(p)^d$.  Finiteness of the ideal class
group gives $\mathfrak a^{N_0}=(\gamma)$ for some $N_0\ge1$ and
$\gamma\in\mathcal O_K$.  Since
$(\gamma\bar\gamma)=(p)^{dN_0}$, the quotient
\[
  \varepsilon=\frac{\gamma\bar\gamma}{p^{dN_0}}
\]
belongs to $\mathcal O_{K^+}^{\times}$.  It is totally positive, because
for every real embedding $\tau:K^+\hookrightarrow\mathbb R$ and an extension
$\widetilde\tau:K\hookrightarrow\C$ one has
$\tau(\gamma\bar\gamma)=|\widetilde\tau(\gamma)|^2>0$.
Hence $\pi=\gamma^2/\varepsilon$ is an algebraic integer.  Put
\[
                  \pi=\gamma^2/\varepsilon,\qquad q=p^{2dN_0}.
\]
Then $\pi\bar\pi=q$; since $K$ is CM, every complex conjugate of
$\pi$ has absolute value $q^{1/2}$.  Thus $\pi$ is a $q$-Weil number.
Since $(\pi)=\mathfrak a^{2N_0}$, its valuations at the four primes,
normalized by $v(p)=1$, are
$(dN_0,dN_0,2N_0,2N_0(d-1))$.  As $d\ge4$, the three numbers $dN_0$,
$2N_0$ and $2N_0(d-1)$ are distinct, so the value $2N_0$ occurs exactly
once; the same holds for every power $\pi^m$, whose valuations are
multiplied by $m$.  By \cref{lem:pair-valuations}, $\Q(\pi^m)=K$ for
all $m\ge1$.

As all four local degrees are one, \eqref{eq:HT-invariant-convention}
gives the stated invariants.  Their least common denominator is $d$,
and $[K:\Q]=4$, so Honda--Tate gives a division algebra of degree $d$
and a variety of dimension $2d$.  The base-change criterion gives
geometric simplicity and equality of the geometric and rational
endomorphism algebras.  All Newton slopes belong to
$\{1/2,1/d,1-1/d\}$, so $B_d$ has $p$-rank zero.
\end{proof}

At the first conjugate pair $(\mathfrak p_1,\bar{\mathfrak p}_1)$ one
has $D_{d,\mathfrak p_1}\simeq M_{d/2}(\Delta_2)$, where
$\Delta_2/\Q_p$ is quaternion division.  Its division index is $2$,
although its degree as a central simple algebra is $d$.  This pair
carries the deformation.  At the second pair
$(\mathfrak p_2,\bar{\mathfrak p}_2)$ the local division index is $d$;
this pair fixes the global degree and remains constant in the family.
Choose a polarization of $B_d$ defined over $\F_q$.  Its Rosati
involution sends Frobenius to Verschiebung, hence
$\pi^\dagger=q/\pi=\bar\pi$.  Since $K=\Q(\pi)$, the Rosati
involution restricts to CM conjugation on $K$ and exchanges the two
local factors of each conjugate pair.

\section{Polarized globalization}\label{sec:globalization}

The following construction inserts one local deformation, together with
its dual, into an abelian scheme.  It applies to the quadratic and
quartic CM seeds, and to the sextic CM seeds of \cref{sec:complements}.

We use the convention of Rapoport and Zink for Cartier duality and
biduality \cite[\S3.20]{RZ}: a quasi-polarization of a $p$-divisible group
is a quasi-isogeny $\lambda:X\to X^\vee$ with
$\lambda^\vee=-\lambda$ under the chosen identification
$X^{\vee\vee}\simeq X$.  With
$A[p^\infty]^\vee\simeq A^\vee[p^\infty]$, an abelian polarization
induces such a map; the sign records the alternating Weil pairing.
All dual maps below use this convention.

\begin{proposition}
\label{prop:quartic-globalization-complete}
Let \((A_0,\lambda_0)/k\) be a polarized abelian variety over an
algebraically closed field of characteristic \(p\).  Let \(D\) be a
central division algebra over a CM field \(K\) of degree \(2e_0\),
equipped with a unital \(\Q\)-algebra homomorphism
\(\iota_0:D\to\Endo(A_0)\).  This homomorphism is injective, and we
identify \(D\) with its image.  No simplicity assumption is imposed on
\(A_0\).  Assume that the Rosati involution stabilizes \(D\) and
restricts to CM conjugation on \(K\), and that \(p\) splits completely
in \(K\).
Write the rational central decomposition as
\[
 K\otimes\Q_p=\prod_{i=1}^{e_0}
 (K_{\mathfrak p_i}\times K_{\bar{\mathfrak p}_i})\simeq\Q_p^{2e_0},
 \qquad
 A_0[p^\infty]\otimes\Q_p
 =\bigoplus_{i=1}^{e_0}(N_i\oplus N_{\bar i}).
\]
Suppose that a \(p\)-divisible group \(Y/k[[t]]\) carries an action of a
full order \(\mathcal O_1\subset D_{\mathfrak p_1}\) and that there is a
\(D_{\mathfrak p_1}\)-linear isomorphism \(Y_0\otimes\Q_p\simeq N_1\).

Then there are a \(p\)-power isogeny \(\psi:A_0\to A'_0\), a
polarization \(\lambda'_0\) with
\(\psi^\vee\lambda'_0\psi=p^{2m}\lambda_0\) for some \(m\ge0\), a
Rosati-stable full order \(\Lambda\subset D\), and a polarized
\(\Lambda\)-linear abelian scheme \((\mathcal A,\lambda)/k[[t]]\)
with special fibre \((A'_0,\lambda'_0)\).  The
\(\mathfrak p_1\)-factor of \(\mathcal A[p^\infty]\) is \(Y\), its
\(\bar{\mathfrak p}_1\)-factor is \(Y^\vee\) with the Rosati-conjugate
action, and all remaining local factors are constant.  The transported \(D\)-action on \(A'_0\) has the original
Rosati involution, and \(D=\Lambda\otimes\Q\) acts faithfully on
every fibre of \(\mathcal A\).
\end{proposition}

\begin{proof}
\smallskip\noindent\emph{Step 1. The local product.}
Let \(e_i,e_{\bar i}\) be the central idempotents of the displayed
decomposition.  Since \(e_i^\dagger=e_{\bar i}\), the identity
\(d^\vee\lambda_0=\lambda_0d^\dagger\), extended to
\(D\otimes\Q_p\), pairs \(N_i\) with \(N_{\bar i}\) and gives
\(N_{\bar i}\simeq N_i^\vee\).  Under these identifications the
polarization is the cross-duality quasi-polarization on each pair.
Put \(Z_1=Y_0\).  For \(i>1\), choose a full order
\(\mathcal O_i\subset D_{\mathfrak p_i}\) and an \(F,V\)-stable
Dieudonn\'e lattice \(M_i\subset\mathbb D(N_i)\).  Replacing it by
\(\mathcal O_iM_i\), the sum of its translates under a finite
\(\Z_p\)-basis of \(\mathcal O_i\), gives a full
\(\mathcal O_i,F,V\)-stable lattice: the order is closed under
multiplication and its elements commute with \(F,V\).
Let \(Z_i\) be the corresponding \(p\)-divisible group.
Give \(Z_{\bar i}=Z_i^\vee\) the action of
\(\mathcal O_{\bar i}=\mathcal O_i^\dagger\) by duality, and set
\(Z=\prod_{i=1}^{e_0}(Z_i\times Z_i^\vee)\).  On each pair use the
standard cross-duality principal quasi-polarization
\(\lambda_{Z,i}(x,\varphi)=(\varphi,-x)\), with the covariant
Cartier-duality convention, and denote their product by \(\lambda_Z\).

\smallskip\noindent\emph{Step 2. The integral special fibre.}
Choose \(D_{\mathfrak p_i}\)-linear quasi-isogenies
\(u_i:N_i\to Z_i\otimes\Q_p\) and put
\(u_{\bar i}=(u_i^\vee)^{-1}\).  The identities
\(u_i^\vee u_{\bar i}=1\) and \(u_{\bar i}^\vee u_i=1\) show,
componentwise on the cross-duality maps, that their product \(u\)
satisfies \(u^\vee\lambda_Zu=\lambda_0\).
Choose \(m\) so that \(\rho=p^m u\) is an integral isogeny.
Its kernel \(H\) is a finite subgroup of \(A_0[p^N]\) for some
\(N\), so the quotient \(\psi:A_0\to A'_0=A_0/H\) is a
\(p\)-power isogeny.  Since \(H\subset A_0[p^\infty]\), the induced
quotient on \(p\)-divisible groups is \(A_0[p^\infty]/H\); by the
construction of \(\rho\), this quotient is \(Z\).  Equivalently,
\(\rho\) factors through \(\psi[p^\infty]\) and induces an isomorphism
\[
 \bar\rho:A'_0[p^\infty]\xrightarrow{\sim} Z.
\]
We use \(\bar\rho\) to identify these two \(p\)-divisible groups.  Define
\[
 \lambda'_0=(\psi^\vee)^{-1}p^{2m}\lambda_0\psi^{-1},
 \qquad \rho^\vee\lambda_Z\rho=p^{2m}\lambda_0.
\]
Some \(f=p^a\lambda'_0\) is integral, because \(\psi\) has
\(p\)-power degree.  Its restriction to the \(p\)-divisible group is
\(f[p^\infty]=p^a\lambda_Z\); this equality holds integrally since
Hom groups of \(p\)-divisible groups have no \(p\)-torsion.
Thus \(f\) kills \(A'_0[p^a]\) scheme-theoretically and factors
through \([p^a]:A'_0\to A'_0\).  Cancellation in the rational Hom
space shows that \(\lambda'_0\) itself is integral.

The rational formula makes \(\lambda'_0\) symmetric.  By
\cite[\S20, Theorem~2]{Mumford}, \(2\lambda'_0=\phi_L\) for a line
bundle \(L\).  Since \(\phi_L=\lambda'_0\circ[2]\), the group scheme
\(A'_0[2]\) is contained in \(\ker\phi_L\) in every characteristic;
hence \cite[\S23, Theorem~3, p.~231]{Mumford} gives \(L\simeq M^{\otimes2}\),
also in characteristic \(2\).  Cancelling \(2\) in the torsion-free
Hom group gives \(\lambda'_0=\phi_M\).
Choose an ample \(L_0\) inducing \(\lambda_0\).  Since
\(\phi_{\psi^*M}=p^{2m}\lambda_0=\phi_{L_0^{\otimes p^{2m}}}\),
the line bundle \(\psi^*M\) is a degree-zero twist of an ample bundle,
and hence ample.  Ampleness descends under the finite surjective map
\(\psi\), so \(M\) is ample and \(\lambda'_0\) is a polarization.
It induces exactly \(\lambda_Z\) on \(Z\).  Transporting the
\(D\)-action through \(\psi\), its Rosati involution satisfies
\[
 (\psi d\psi^{-1})^{\dagger'}=\psi d^\dagger \psi^{-1}\qquad(d\in D).
\]
We identify the transported algebra and involution with \(D,\dagger\).

\smallskip\noindent\emph{Step 3. A Rosati-stable order.}
Set
\[
 \begin{split}
 \Lambda_0=\{d\in D\cap\End(A'_0):\;&
 d_{\mathfrak p_i}\in\mathcal O_i,\quad
 d_{\bar{\mathfrak p}_i}\in\mathcal O_i^\dagger
 \text{ for }1\le i\le e_0\},\\
 \Lambda&=\Lambda_0\cap\Lambda_0^\dagger.
 \end{split}
\]
The intersection \(D\cap\End(A'_0)\) is a full \(\Z\)-order in
\(D\), and the local conditions impose finite-index lattice
conditions only at \(p\).  Thus \(\Lambda_0\), and then \(\Lambda\),
are full orders.  Under the identification \(\bar\rho:A'_0[p^\infty]\simeq Z\), the
projections of \(\Lambda\otimes\Z_p\) act through the corresponding
integral local actions on the factors of \(Z\).  Their actions on
conjugate factors are adjoint under \(\lambda_Z\).

\smallskip\noindent\emph{Step 4. Serre--Tate theory and algebraization.}
Put \(R=k[[t]]\) and form the deformation
\[
 \mathcal Z=Y\times Y^\vee\times
 \prod_{i=2}^{e_0}\bigl((Z_i\times Z_i^\vee)\times_kR\bigr),
\]
with its product cross-duality quasi-polarization
\(\lambda_{\mathcal Z}\).  The local actions give a \(\Lambda\)-action
\(\iota_{\mathcal Z}\) satisfying
\[
 \iota_{\mathcal Z}(a)^\vee\lambda_{\mathcal Z}
 =\lambda_{\mathcal Z}\iota_{\mathcal Z}(a^\dagger)
 \qquad(a\in\Lambda).
\]
For \(R_j=R/(t^{j+1})\), Serre--Tate theory
\cite[Theorem~1.2.1]{KatzST} gives an abelian scheme
\(\mathcal A_j/R_j\), uniquely up to unique isomorphism, whose special
fibre is \(A'_0\) and whose \(p\)-divisible group is
\(\mathcal Z_{R_j}\).  Katz's theorem gives an equivalence of categories,
so it is fully faithful.  Together with compatibility with duality, this
lifts the \(\Lambda\)-action and the pair
\((\lambda'_0,\lambda_{\mathcal Z}|_{R_j})\) uniquely to a symmetric
homomorphism
\[
 \lambda_j:\mathcal A_j\longrightarrow\mathcal A_j^\vee.
\]
The Rosati identities lift as well.  Uniqueness makes these objects and
morphisms compatible as \(j\) varies.

Let \(\mathcal P_j\) be the rigidified Poincar\'e bundle on
\(\mathcal A_j\times_{R_j}\mathcal A_j^\vee\), and put
\[
 \mathcal L_j=(1,\lambda_j)^*\mathcal P_j.
\]
The Poincar\'e identity gives \(\phi_{\mathcal L_j}=2\lambda_j\).
On the special fibre this is \(2\lambda'_0\), a polarization; hence
\(\mathcal L_j|_{A'_0}\) is ample.  Ampleness is unchanged by a
nilpotent thickening, so \(\mathcal L_j\) is relatively ample over
\(R_j\).  The bundles \(\mathcal L_j\) are compatible under base
change.  This construction, including characteristic \(2\), is the
usual Poincar\'e construction over nilpotent bases; see
\cite[Ch.~I, \S1]{FaltingsChai}.

The compatible system \((\mathcal A_j,\mathcal L_j)_j\) is therefore a
proper formal scheme over \(\operatorname{Spf} R\) endowed with an ample invertible
sheaf.  Grothendieck algebraization \cite[5.4.5]{EGAIII} produces a
projective scheme \(\mathcal A/R\) with relatively ample line bundle
\(\mathcal L\) whose formal completion is the given system.  Formal
GAGA \cite[5.1.6 and 5.4.1]{EGAIII}, applied to the graphs of the
compatible formal morphisms, algebraizes the multiplication, identity,
inverse, and a finite \(\Z\)-basis of the \(\Lambda\)-action.  All
group and ring relations hold because they hold modulo every
\(t^{j+1}\).  Thus \(\mathcal A\) is a proper group scheme over \(R\).

We now verify that it is an abelian scheme.  At every point of the
special fibre, the completed local morphism is the one defined by the
formal abelian scheme \((\mathcal A_j)_j\), hence is formally smooth.
By the formal criterion for smoothness, \(\mathcal A\to\Spec R\) is
smooth along the special fibre.  Its nonsmooth locus is closed, since
\(\mathcal A/R\) is of finite presentation.  If it were nonempty, its
image under the proper morphism \(\mathcal A\to\Spec R\) would be a
nonempty closed subset disjoint from the closed point.  No such subset
exists in \(\Spec k[[t]]\).  Hence \(\mathcal A/R\) is smooth
everywhere.  For a proper smooth morphism the number of geometric connected
components of the fibres is locally constant
\cite[Corollaire~7.8.10(i)]{EGAIII2}.  The special
fibre is the connected abelian variety \(A'_0\), and \(\Spec R\) is
connected; hence every geometric fibre is connected.  Therefore
\(\mathcal A/R\) is an abelian scheme.

Now \(\mathcal A^\vee\) exists and commutes with base change.  The
compatible maps \(\lambda_j\) algebraize, again by formal GAGA, to a
homomorphism
\[
 \lambda:\mathcal A\longrightarrow\mathcal A^\vee.
\]
It is symmetric, and the algebraized line bundle satisfies
\(\phi_{\mathcal L}=2\lambda\).  We claim that \(\lambda\) is a
polarization.  On any geometric fibre \(A_s\), the symmetric map
\(\lambda_s\) satisfies \(2\lambda_s=\phi_{\mathcal L_s}\), with
\(\mathcal L_s\) ample.  As in Step~2,
\(A_s[2]\subseteq\ker\phi_{\mathcal L_s}\); hence
\cite[\S23, Theorem~3, p.~231]{Mumford} gives
\(\mathcal L_s\simeq N_s^{\otimes2}\) for a line bundle \(N_s\).
Cancelling \(2\) in the torsion-free group
\(\operatorname{Hom}(A_s,A_s^\vee)\) yields \(\lambda_s=\phi_{N_s}\), and
\(N_s\) is ample because \(\mathcal L_s\) is ample.  Thus every
geometric fibre of \(\lambda\) is a polarization, so \(\lambda\) is a
polarization of the abelian scheme.

It remains to identify the \(p\)-divisible group.  For every \(n\ge1\)
and every \(j\), the Serre--Tate construction gives compatible
isomorphisms
\[
 \mathcal A[p^n]_{R_j}\simeq\mathcal Z[p^n]_{R_j}.
\]
Both sides are finite flat over the complete ring \(R\).  Their finite
coordinate algebras are \(t\)-adically complete, so the compatible
isomorphisms modulo \(t^{j+1}\) pass uniquely to the inverse limit and
give
\[
 \mathcal A[p^n]\simeq\mathcal Z[p^n].
\]
These isomorphisms are compatible in \(n\), and hence
\(\mathcal A[p^\infty]\simeq\mathcal Z\).  The Rosati identities hold
after algebraization because they hold modulo every \(t^{j+1}\).
Finally \(\Lambda\otimes\Q=D\).  On every geometric fibre the induced
map \(D\to\Endo(\mathcal A_s)\) is unital; its kernel is a two-sided
ideal of the division algebra \(D\), and is therefore zero.  Thus the
\(D\)-action is faithful on every fibre.
\end{proof}

\section{Realization of the ratios}\label{sec:realization}

We first apply the quaternion window to obtain the uniform
$p$-rank-zero result.  The degree-$s$ window then gives a smaller
quadratic-CM realization with positive $p$-rank.

\begin{theorem}\label{thm:uniform-quartic}
Let $p$ be a prime, let $d\ge4$ be an even integer, and let $s>r\ge3$
be coprime integers; put $n=r+s$.
There exist an algebraically closed field $\Omega$ of characteristic
$p$ and a simple abelian variety $A/\Omega$ with
\[
 \dim A=2nd,\qquad f(A)=0,\qquad m(A)=\frac nd.
\]
Its geometric endomorphism algebra is the algebra $D_d$ of
\cref{prop:quartic-seed-complete}, a central division algebra of degree $d$
over the quartic CM field $K$ of \cref{lem:quartic-cm-split-complete}.
Moreover, $A$ is the geometric generic fibre of an abelian scheme over
$\overline{\F}_p[[t]]$ all of whose geometric fibres have $p$-rank zero.
\end{theorem}

\begin{proof}
\smallskip\noindent\emph{Step 1. The seed and the local deformation.}
The hypotheses on $d$ are those of \cref{prop:quartic-seed-complete},
and those on $r,s$ include those of
\cref{prop:two-level-window-complete}; the two choices are independent.
Choose the seed $B_d/\F_q$ and write $D=D_d$ and $K=Z(D)$.
Work over $k=\overline{\F}_p$ and put $R=k[[t]]$.
Equip $B_d$ with a polarization defined over $\F_q$ and
$A_0=B_d^n$ with its product polarization.
The diagonal $D$-action is Rosati-stable; since
$\pi^\dagger=q/\pi=\bar\pi$ and $K=\Q(\pi)$, the Rosati involution
induces CM conjugation on $K$.

Let $X=X_2/R$ be the quaternion group of
\cref{prop:two-level-window-complete}, and put $Y=X^{d/2}$.
The standard matrix action makes $Y$ linear over
$M_{d/2}(\mathcal O_2)\subset
D_{\mathfrak p_1}\simeq M_{d/2}(\Delta_2)$.
It has height $nd$, dimension $nd/2$, and isoclinic special fibre
of slope $1/2$.
The $\mathfrak p_1$-factor of $A_0[p^\infty]$ has the same height
and special slope.  Fix an underlying isocrystal isomorphism.
The two actions of $D_{\mathfrak p_1}$ are embeddings into the
central simple algebra $M_{nd/2}(\Delta_2)$, agreeing on $\Q_p$;
Skolem--Noether makes them conjugate.  Adjusting the isomorphism
therefore identifies $Y_0$ with that central factor
$D_{\mathfrak p_1}$-linearly, up to isogeny.

\smallskip\noindent\emph{Step 2. Globalization and the endomorphism algebra.}
By \cref{prop:quartic-globalization-complete}, there is a polarized
abelian scheme $\mathcal A/R$ carrying a rational $D$-action,
with special fibre $D$-linearly $p$-power isogenous to $B_d^n$,
active factors $Y,Y^\vee$, and the other conjugate pair constant.
The reduced geometric generic slopes of $X$ are
$(1-1/r)^{[r]}$ and $(1+1/s)^{[s]}$.
Their multiplicities equal their denominators in lowest terms.
Reduced centralizers and the matrix action in \cref{prop:rank2}
give
\[
 \Endo_{D_{\mathfrak p_1}}(Y_{\bar\eta})
 \simeq\mathcal D_r\times\mathcal D_s,
 \qquad [\mathcal D_r:\Q_p]=r^2,\quad
 [\mathcal D_s:\Q_p]=s^2.
\]
These factors act on the full-slope components of heights $dr,ds$
and slopes $1/2-1/r,1/2+1/s$, respectively.  Both slopes differ from
the special slope $1/2$, and $\gcd(r,s)=1$, so all hypotheses of
\cref{prop:commutant-criterion} hold and
$\Endo(\mathcal A_{\bar\eta})=D$.  Concretely, an idempotent selecting
the first component would, after the isogeny of
\cref{lem:spread-idempotent}, give a summand with generic dimension
$dr(1/2+1/r)=dr/2+d$ and special dimension $dr/2$; for the second
component these dimensions would be $ds(1/2-1/s)=ds/2-d$ and $ds/2$.
In particular, $A=\mathcal A_{\bar\eta}$ is simple.

\smallskip\noindent\emph{Step 3. The $p$-rank and the ratio.}
The two active factors have full slopes
$1/2\pm1/r$ and $1/2\pm1/s$; the constant pair has slopes
$1/d$ and $1-1/d$.  All lie strictly between zero and one because
$s>r\ge3$ and $d\ge4$.  Consequently $f(A)=0$.
The special fibre is isogenous to $B_d^n$ and also has $p$-rank
zero.  Thus both the special and geometric generic fibres have $p$-rank
zero; since $\Spec k[[t]]$ has only its closed and generic points, every
geometric fibre does.  Finally, constancy of dimension and the degree of
$D$ give
\[
 \dim A=n\dim B_d=2nd,
 \qquad
 m(A)=\frac{4nd}{4d^2}=\frac nd.
\]
Take $\Omega=\overline{k((t))}$.
\end{proof}

To realize a prescribed ratio we need only the following elementary
fact.

\begin{lemma}\label{lem:coprime-split}
An integer $N$ can be written as $N=r+s$ with coprime integers
$s>r\ge3$ if and only if $N\ge7$.
\end{lemma}

\begin{proof}
Such a decomposition forces $N\ge3+4$.  Conversely, for odd $N\ge7$
take $(r,s)=((N-1)/2,(N+1)/2)$, which are consecutive.  For
$N\equiv0\pmod4$ and $N\ge8$ take $(N/2-1,N/2+1)$: both are odd and
they differ by $2$.  For $N\equiv2\pmod4$ and $N\ge10$ take
$(N/2-2,N/2+2)$: both are odd and they differ by $4$.
\end{proof}

\begin{corollary}\label{cor:ratio-dimension}
Let $p$ be a prime, let $a\ge2$ and $b\ge1$ be integers, and let
$c\ge1$ be an integer such that $cb$ is even, $cb\ge4$ and $ca\ge7$.
There exists a simple abelian variety $A$ over an algebraically closed
field of characteristic $p$ with
\[
 \dim A=2c^2ab,\qquad f(A)=0,\qquad m(A)=\frac ab.
\]
We can choose $c=4$, which gives $\dim A=32ab$.  If $a\ge4$ and
$b\ge2$, we can choose $c=2$, giving $\dim A=8ab$.  If $a\ge7$ and
$b\ge4$ is even, we can choose $c=1$, giving $\dim A=2ab$.
\end{corollary}

\begin{proof}
Put $d=cb$ and $n=ca$, write $n=r+s$ by \cref{lem:coprime-split},
and apply \cref{thm:uniform-quartic}; then $\dim A=2nd=2c^2ab$ and
$m(A)=n/d=a/b$.  Each of the three displayed choices satisfies the
hypotheses.  For $c=4$, \cref{lem:coprime-split} gives
$(r,s)=(2a-1,2a+1)$.
\end{proof}

\begin{proof}[Proof of \cref{thm:main}]
If $\mu=1/b$ with $b\ge3$, the seed $B_{p,b}$ of \cref{prop:seeds},
base-changed to $\overline{\F}_p$, is simple of $p$-rank zero with
$m=1/b$.  If $\mu=1/2$, take a supersingular elliptic curve: its
geometric endomorphism algebra is the quaternion algebra over $\Q$
ramified at $p$ and $\infty$, so $m=2/4=1/2$, and its $p$-rank is zero.
Every other $\mu\in\Q_{>0}$ has a representation $a/b$ with $a\ge2$
and $b\ge1$, using $1=2/2$; apply \cref{cor:ratio-dimension} with
$c=4$.
\end{proof}

\begin{remark}
\label{rem:classical-prank-zero}
The unit fractions $1/b$ with $b\ge2$ are realized with $p$-rank zero by
a supersingular elliptic curve and by the seeds of \cref{prop:seeds}.
The ratio $1$ is realized with $p$-rank zero by the construction of
Lenstra and Oort with $t=1$, $(n_1,m_1)=(2,1)$ and $h=0$: it gives a
simple abelian threefold over $\overline{\F}_p$ with formal isogeny type
$G_{2,1}+G_{1,2}$, whose geometric endomorphism algebra is a CM field of
degree six by \cite[(2.13) and Lemma~2.14]{LenstraOort}, applied over
every finite extension of the field of definition.  For $g\ge3$, the
symmetric Newton polygon
$(1/3)^{[3]}\cup(1/2)^{[2g-6]}\cup(2/3)^{[3]}$ has $p$-rank zero and is
not supersingular.  For $N\ge3$ prime to $p$, its Newton stratum in
$\mathcal A_{g,1,N}$ is nonempty by Manin's conjecture \cite{OortNP}
(see also \cite[p.~314]{OortQ}), geometrically irreducible by
\cite[Theorem~3.1]{ChaiOortMonodromy}, and stable under prime-to-$p$
Hecke correspondences.  Applying \cite[Corollary~4.2]{ChaiMonodromy} to
its smooth locus shows that the abelian variety over its geometric
generic point has endomorphism ring $\Z$; it is simple with $f=0$ and
$m=2g$.  The realization of the endomorphism ring $\Z$ for a prescribed
formal isogeny type in \cite[Corollary~1.4]{OortVanDerPut} assumes that
the type has positive $p$-rank.
\end{remark}

\subsection{Choice of parameters}\label{sec:parameters}

\Cref{tab:choices} lists one realization of $p$-rank zero for each
ratio.  For $\mu=1$ it lists the Lenstra--Oort threefold, whereas the
proof of \cref{thm:main} writes $1=2/2$ and uses the general
construction.  In the last row, $K$ is the quartic CM field
of \cref{lem:quartic-cm-split-complete}, the algebra $D=D_d$ has
invariants $(1/2,1/2,1/d,1-1/d)$ at the four primes above $p$, the
integer $n$ is the exponent of the seed, and $g=\dim A$.

\begin{table}[ht]
\centering
\renewcommand{\arraystretch}{1.25}
\begin{tabular}{l|l|c|c}
Target $\mu$ & Construction and parameters & $[\Endo(A):\Q]$ & $g$\\
\hline
$1/2$ & supersingular elliptic curve & $4$ & $1$\\
$1/b$, $b\ge3$ & quadratic seed $B_{p,b}$ of \cref{prop:seeds} & $2b^2$ & $b$\\
$1$ & Lenstra--Oort threefold (\cref{rem:classical-prank-zero}) & $6$ & $3$\\
$a/b$, $a\ge2$ & $n=ca$, $d=cb$, with $d\ge4$ even and $n\ge7$ & $4c^2b^2$ & $2c^2ab$
\end{tabular}
\caption{Choices of $p$-rank zero, valid for every prime $p$.}
\label{tab:choices}
\end{table}

For the last row, write $\mu=a/b$ in lowest terms and put
\[
 L=\max\Bigl\{1,\Bigl\lceil\frac4b\Bigr\rceil,
       \Bigl\lceil\frac7a\Bigr\rceil\Bigr\},\qquad
 c_{\min}=
 \begin{cases}
 L,& b\text{ even},\\
 2\lceil L/2\rceil,& b\text{ odd}.
 \end{cases}
\]
This is the smallest multiplier satisfying the hypotheses of
\cref{cor:ratio-dimension}.  After taking $c=c_{\min}$, the explicit
rule in \cref{lem:coprime-split} gives $n=r+s$ with $\gcd(r,s)=1$ and
$s>r\ge3$.  \Cref{tab:examples} lists some examples.  The dimensions are
those produced by the construction; they are not claimed to be minimal.

\begin{table}[ht]
\centering
\renewcommand{\arraystretch}{1.2}
\begin{tabular}{c|cccccc}
$\mu$ & $c_{\min}$ & $n$ & $d$ & $r$ & $s$ & $\dim A$\\
\hline
$2/3$ & $4$ & $8$ & $12$ & $3$ & $5$ & $192$\\
$3/2$ & $3$ & $9$ & $6$ & $4$ & $5$ & $108$\\
$4/3$ & $2$ & $8$ & $6$ & $3$ & $5$ & $96$\\
$5/4$ & $2$ & $10$ & $8$ & $3$ & $7$ & $160$\\
$7/4$ & $1$ & $7$ & $4$ & $3$ & $4$ & $56$\\
$2$ & $4$ & $8$ & $4$ & $3$ & $5$ & $64$
\end{tabular}
\caption{Examples of the quartic construction.}
\label{tab:examples}
\end{table}

Every row of \cref{tab:examples} has $f(A)=0$ and $\Endo(A)=D_d$.  The
table minimizes $c$ only within \cref{cor:ratio-dimension}; other
constructions may have smaller dimension.  For example, the quadratic
construction of \cref{prop:quadratic-realization} realizes $5/3$ in
dimension $15$ with $p$-rank $9$.

For a reduced fraction $a/b$ with $a\ge7$ and even $b\ge4$, taking
$c=1$ in \cref{cor:ratio-dimension} gives a family of $p$-rank zero
whose special fibre is $D_b$-linearly $p$-power isogenous to $B_b^a$,
where $m(B_b)=1/b$.  This gives the quartic-CM range related to
Question~22.3, with the polarization and endomorphism order modified
as in \cref{prop:quartic-globalization-complete}.

\begin{remark}\label{rem:field-of-definition}
In \cref{thm:main} one may take $k$ to be an algebraic closure of a
finitely generated extension $F$ of $\F_p$.  Indeed, by standard limit
arguments \cite[\S8]{EGAIV3}, the generic fibre of the abelian scheme
in \cref{thm:uniform-quartic} is the base change of an abelian variety
$A_F$ over a finitely generated subfield $F\subset\overline{\F}_p((t))$.
An endomorphism of $A_F\otimes\Omega$ is defined over a finitely
generated $\overline F$-subalgebra of $\Omega$, and therefore gives a
family of endomorphisms of $A_{\overline F}$ parametrized by a connected
$\overline F$-variety; by the rigidity lemma \cite[Proposition~6.1]{GIT} such a family
is constant.  Hence $\Endo(A_{\overline F})=\Endo(A_F\otimes\Omega)$,
and geometric simplicity, $m$ and the $p$-rank are unchanged under this
extension of algebraically closed base fields.  When the reduced
numerator of $\mu$ exceeds one, $F$ is necessarily transcendental over
$\F_p$ by the Honda--Tate discussion in the introduction.  We do not
know whether, for every ratio, one can take
$\operatorname{trdeg}(F/\F_p)=1$.
\end{remark}

\begin{remark}\label{rem:dimension}
\Cref{cor:ratio-dimension} realizes $a/b$ with $p$-rank zero in
dimension $2c^2ab$, and \cref{prop:quadratic-realization} below gives
dimension $ab$ in its range, with positive $p$-rank.  These dimensions
are effective outputs of the constructions, not optimal values.  We do
not know the smallest dimension of a simple abelian variety $A$ with
$f(A)=0$ and $m(A)=a/b$.
\end{remark}

\begin{example}\label{ex:mu-two-thirds}
Take $a=2$ and $b=3$.  Here $c=4$ is the smallest value satisfying the
hypotheses of \cref{cor:ratio-dimension}, so $d=12$, $n=8$, $r=3$ and
$s=5$.
The geometric generic fibre of the quaternion window $X=X_2$ has
reduced slopes $(2/3)^{[3]}\cup(6/5)^{[5]}$ and full covariant slopes
$(1/6)^{[6]}\cup(7/10)^{[10]}$.
Forming $Y=X^6$, adjoining its dual, and including the constant
pair gives the full polygon of $A$:
\[
 (1/12)^{[96]}\cup(1/6)^{[36]}\cup(3/10)^{[60]}
 \cup(7/10)^{[60]}\cup(5/6)^{[36]}\cup(11/12)^{[96]}.
\]
There is no slope zero or one, so $f(A)=0$.
The seed $B_{12}$ has dimension $24$, and $A$ has dimension $192$.
Its geometric endomorphism algebra has $\Q$-dimension
$4\cdot12^2=576$, giving $m(A)=2\cdot192/576=2/3$.

For this ratio the integers $n,d$ of \cref{thm:uniform-quartic} satisfy
$3n=2d$.  Since $d$ is even, write $n=2c$ and $d=3c$ with $c$ even.
The condition $n\ge7$ then forces $c\ge4$, which explains the value
$192$ within this construction.  This is an upper bound, not a
minimality result (\cref{rem:dimension}).
\end{example}

\subsection{The smaller quadratic-CM construction}

\begin{proposition}
\label{prop:quadratic-realization}
Let $p$ be a prime and let $s>r\ge2$ be coprime integers, with $n=r+s$.
There exist an algebraically closed field $\Omega$ of characteristic $p$
and a simple abelian variety $A/\Omega$ such that
\[
 \dim A=ns,\qquad \Endo(A)\simeq D_{p,s},\qquad
 m(A)=\frac ns,\qquad f(A)=s^2,
\]
where $D_{p,s}$ is the division algebra of \cref{prop:seeds}.
\end{proposition}

\begin{proof}
Let $B=B_{p,s}$, $D=D_{p,s}$, and $K=K_p$, and work over an algebraic
closure $k$ of the finite field of definition of $B$.
Let $\pfr\mid p$ be the prime at which the covariant slope of
$B[p^\infty]$ is $1-1/s$; with the convention of \cref{sec:seeds},
$\pfr=(\alpha)$.  The local algebra $D_{\pfr}$ is a central division
algebra of degree $s$ over $\Q_p$.

Let $X_s/k[[t]]$ be the degree-$s$ group of \cref{thm:gen-window}.  Its
special fibre and the $\pfr$-factor of $B^n[p^\infty]$ are isoclinic of
height $ns$ and slope $1-1/s$.  Fix an isomorphism $\theta$ of their
rational Dieudonn\'e modules, and let $E_N$ be their common
quasi-endomorphism algebra, central simple over $\Q_p$ of degree $ns$.
Through $\theta$, both $D_{\pfr}$ and the algebra $\Delta_s$ of
\cref{thm:gen-window} embed unitally into $E_N$.  Let $C_D$ and $C_\Delta$ be the centralizers of $D_{\pfr}$ and
$\Delta_s$ in $E_N$, respectively.  By the double-centralizer theorem,
both are central simple of degree $n$, with Brauer classes
\[
 [C_D]=[E_N]-[D_{\pfr}],\qquad
 [C_\Delta]=[E_N]-[\Delta_s].
\]
Their periods divide their indices, hence divide $n$.  Subtracting the
two relations shows that $n([\Delta_s]-[D_{\pfr}])=0$.  On the other
hand, both $[D_{\pfr}]$ and $[\Delta_s]$ are annihilated by $s$, since
the two algebras have degree $s$.  Thus the difference
$[\Delta_s]-[D_{\pfr}]$ is killed by both $n$ and $s$.  Since
$\gcd(n,s)=\gcd(r,s)=1$, it is zero.  The two algebras are central
division algebras of the same degree $s$, and therefore
$D_{\pfr}\simeq\Delta_s$.
After choosing such an identification, Skolem--Noether conjugates the
two embeddings into $E_N$; replacing $\theta$ by $g\theta$ for a suitable
$g\in E_N^\times$ therefore makes it $D_{\pfr}$-linear.
The product polarization on $B^n$ has Rosati-stable diagonal $D$-action,
with CM conjugation on $K$.  The split-CM globalization
\cref{prop:quartic-globalization-complete} therefore gives an abelian
scheme $\mathcal A/k[[t]]$ with rational $D$-action, special fibre
$D$-linearly isogenous to $B^n$, and local factors $X_s$ and $X_s^\vee$.

The generic reduced slopes of $X_s$ have multiplicities equal to their
coprime denominators $r,s$.  By \cref{prop:rank2}, the geometric generic
centralizer $\Endo_{D_{\pfr}}((X_s)_{\bar\eta})$ is a product of division
algebras of degrees $r,s$.
Put $\lambda=1-1/r-1/s$.  The full generic components of $X_s$ have
heights $rs,s^2$ and slopes $\lambda,1$.  Their generic dimensions are
respectively
\[
 rs(1-\lambda)=r+s=n,\qquad s^2(1-1)=0,
\]
whereas a special-fibre summand of the same height has dimension
$r$ or $s$, because the special slope is $1-1/s$.
These mismatches exclude both possible nontrivial idempotents in the
completed commutant.  The criterion of \cref{prop:commutant-criterion},
using $\gcd(r,s)=1$, gives $\Endo(\mathcal A_{\bar\eta})=D$.
Thus $A=\mathcal A_{\bar\eta}$ is simple.

Constancy of dimension gives $\dim A=n\dim B=ns$, and
$[D:\Q]=2s^2$ gives $m(A)=n/s$.
Adding the Cartier-dual factor, the full generic slopes are
\[
 0^{[s^2]}\cup\lambda^{[rs]}
 \cup(1-\lambda)^{[rs]}\cup1^{[s^2]}.
\]
Since $0<\lambda<1$, the $p$-rank is exactly $s^2$.
\end{proof}

For a reduced fraction $a/b$ with $b+2\le a\le2b-1$, take
$s=b$ and $r=a-b$.  The proposition gives dimension $ab$ and
$p$-rank $b^2$, with special fibre isogenous to $B_{p,b}^a$ and
$m(B_{p,b})=1/b$.  This gives the quadratic-CM range related to
Question~22.3, subject to the changes of special fibre, polarization,
and endomorphism order described in the introduction.  The slope-$1$
component of $X_s$ is the reason why the $p$-rank-zero theorem uses the
quaternion realization instead.

\section{Realizations with positive \texorpdfstring{$p$}{p}-rank}
\label{sec:complements}

The construction of \cref{thm:uniform-quartic} keeps all Newton slopes in the
open interval $(0,1)$.  In this section we show that every positive
rational number is also the endomorphism ratio $m(A)$ of a simple
abelian variety $A$ of positive $p$-rank. For nonintegral ratios, however, maximal \(p\)-rank is impossible, as the following standard observation shows.

\begin{remark}\label{rem:ordinary}
If $A$ is a simple ordinary abelian variety over an algebraically closed
field $k$ of characteristic $p$, then $m(A)$ is a positive integer.
Indeed, over each $W_n(k)$ the canonical lifting of $A$ is the lifting
all of whose Serre--Tate parameters equal $1$ \cite[\S4.5]{KatzST};
hence, by \cite[Theorem~2.1(3),(4)]{KatzST}, every endomorphism and
every polarization of $A$ lifts, uniquely, to it. These lifts are
compatible in $n$, and the polarized formal lift over $W(k)$
algebraizes. This gives a unital action of $D=\Endo(A)$ on an abelian
variety $\tilde A$ of dimension $\dim A$ over the fraction field of
$W(k)$. Descend $\tilde A$ and the action of a finite $\Q$-basis of $D$
to a finitely generated field over $\Q$, and embed that field in $\C$.
Then $H_1(\tilde A_{\C},\Q)$ is a vector space over the division
algebra $D$, so $[D:\Q]\mid2\dim A$ and $m(A)\in\Z_{>0}$.
\end{remark}

There are two sources of an \'etale part.  In
\cref{prop:unit-positive} and \cref{thm:positive-rank} it is constant:
one conjugate pair of primes above $p$ carries invariants $0$, so that
the slopes there are $1$ and $0$.  In the covariant convention of
Section~\ref{sec:strategy}, the slope-$1$ part is the \'etale part, so
this pair contributes positive $p$-rank.  For unit fractions this
already happens over a finite field, with the quartic field $K$ of
\cref{lem:quartic-cm-split-complete}.  For general ratios we replace
the quartic centre of \cref{thm:uniform-quartic} by a sextic CM field:
one conjugate pair carries the deformation, one fixes the global
degree, and the third contributes the \'etale part; the rest of the
argument is unchanged.  In \cref{cor:r2-variant}, by contrast, the
quartic seed of $p$-rank zero is kept and the \'etale part is created
by the deformation itself: taking $r=2$ in the quaternion window gives
a generic slope $0$, whose dual has slope $1$.

\begin{lemma}\label{lem:sextic-field}
For every prime $p$ there exists a sextic CM field $K'$ in which $p$
splits completely.  One can take $K'=K_pK_3$, where $K_p$ is the field
of \cref{prop:seeds} and $K_3$ is a cyclic cubic subfield of a cyclotomic
field.  More precisely, one may choose a prime $\ell\equiv1\pmod3$,
$\ell\ne p$, such that $p$ is a cube modulo $\ell$, and take $K_3$ to be
the cubic subfield of $\Q(\zeta_\ell)$.  For $p=2$, one may take
$\ell=31$.
\end{lemma}

\begin{proof}
By the Chebotarev density theorem, there are infinitely many primes
$\ell\nmid 3p$ that split completely in
\[
L=\Q(\zeta_3,\sqrt[3]{p}).
\]
For such a prime, $\ell\equiv1\pmod3$ and $p$ is a cube modulo $\ell$. 
Put
\[
G=(\Z/\ell\Z)^\times,\qquad H=G^3.
\]
Since $G$ is cyclic and $3\mid |G|$, the subgroup $H$ has index $3$.
Hence
\[
K_3=\Q(\zeta_\ell)^H
\]
is a cyclic cubic field.  Moreover, $-1=(-1)^3\in H$, so complex
conjugation acts trivially on $K_3$.  Thus $K_3$ is totally real.

Since $p$ is a cube modulo $\ell$, its class belongs to $H$.  Therefore
the Frobenius of $p$ in
\[
\Gal(K_3/\Q)\simeq G/H
\]
is trivial, and $p$ splits completely in $K_3$.  By
\cref{prop:seeds}, $p$ also splits completely in the imaginary quadratic
field $K_p$.

The degrees of $K_p$ and $K_3$ are $2$ and $3$, so
$K_p\cap K_3=\Q$.  Therefore
\[
[K_pK_3:\Q]=6.
\]
Since $K_3$ is totally real and $K_p$ is imaginary quadratic,
$K'=K_pK_3$ is a totally imaginary quadratic extension of $K_3$ and
hence a CM field.  The Frobenius of $p$ in $\Gal(K'/\Q)$ restricts
trivially to both $K_p$ and $K_3$, and therefore is trivial.  Thus $p$
splits completely in $K'$.

For $p=2$, the choice $\ell=31$ works, since
$4^3\equiv2\pmod{31}$.
\end{proof}

For a sextic CM field $K'$ in which $p$ splits completely, write
$p\mathcal O_{K'}=\prod_{i=1}^3\mathfrak p_i\bar{\mathfrak p}_i$, with
the bar denoting CM conjugation; this labelling is possible for the same
reason as after \cref{lem:pair-valuations}.

\begin{proposition}
\label{prop:unit-positive}
Let $K$ be the quartic CM field of \cref{lem:quartic-cm-split-complete}.
For every integer $d\ge2$ there exist a finite field $\F_q$ of
characteristic $p$ and a geometrically simple abelian variety
$B^\circ_d/\F_q$ of dimension $2d$ and $p$-rank $d$ whose geometric
endomorphism algebra is a central division algebra of degree $d$ over $K$, with
invariants $(0,0,1/d,1-1/d)$ at
$(\mathfrak p_1,\bar{\mathfrak p}_1,\mathfrak p_2,\bar{\mathfrak p}_2)$.
In particular, $m(B^\circ_d\otimes\overline{\F}_p)=1/d$.
\end{proposition}

\begin{proof}
Set $\mathfrak a=\bar{\mathfrak p}_1^{\,d}\mathfrak p_2
\bar{\mathfrak p}_2^{\,d-1}$, so $\mathfrak a\bar{\mathfrak a}=(p)^d$.
As in the proof of \cref{prop:quartic-seed-complete}, choose $N_0\ge1$,
$\gamma\in\mathcal O_K$ and a totally positive unit $\varepsilon$ with
$\mathfrak a^{N_0}=(\gamma)$ and $\gamma\bar\gamma=\varepsilon p^{dN_0}$,
and put $\pi=\gamma^2/\varepsilon$ and $q=p^{2dN_0}$.  The valuations of
$\pi$ at the four primes are $(0,2dN_0,2N_0,2N_0(d-1))$.  The value $0$
occurs exactly once, also for every power $\pi^m$, so
$\Q(\pi^m)=K$ for all $m\ge1$ by \cref{lem:pair-valuations}.
By \eqref{eq:HT-invariant-convention}, the invariants are
$(0,1,1/d,1-1/d)$, that is, $(0,0,1/d,1-1/d)$ in $\Q/\Z$; their least
common denominator is $d$.  Honda--Tate theory gives a division algebra
of degree $d$ over $K$ and a variety of dimension $2d$, and the
base-change criterion of \cref{sec:seeds} gives geometric simplicity and
equality of the geometric and rational endomorphism algebras.  The
covariant slopes at the four primes are $1$, $0$, $1-1/d$ and $1/d$,
each with multiplicity $d$, so the $p$-rank is $d$.  Finally,
$m=2\cdot2d/(4d^2)=1/d$.
\end{proof}

\begin{proposition}\label{prop:sextic-seed}
Let $K'$ be a sextic CM field in which $p$ splits completely, with the
primes above $p$ labelled as above.  For every even integer $d\ge2$ there exist a finite field $\F_q$ of
characteristic $p$ and a geometrically simple abelian variety
$B'_d/\F_q$ of dimension $3d$ and $p$-rank $d$ whose geometric
endomorphism algebra $D'_d$ is a central division algebra of degree $d$ over $K'$,
with invariants
\[
 \Bigl(\frac12,\frac12,\frac1d,1-\frac1d,0,0\Bigr)
 \quad\text{at}\quad
 (\mathfrak p_1,\bar{\mathfrak p}_1,\mathfrak p_2,\bar{\mathfrak p}_2,
  \mathfrak p_3,\bar{\mathfrak p}_3),
 \qquad [D'_d:\Q]=6d^2.
\]
\end{proposition}

\begin{proof}
Set
$\mathfrak a=\mathfrak p_1^{d/2}\bar{\mathfrak p}_1^{d/2}
\mathfrak p_2\bar{\mathfrak p}_2^{d-1}\bar{\mathfrak p}_3^{\,d}$, so
$\mathfrak a\bar{\mathfrak a}=(p)^d$.  Construct $N_0$, $\gamma$,
$\varepsilon$, $\pi=\gamma^2/\varepsilon$ and $q=p^{2dN_0}$ as in the
proof of \cref{prop:quartic-seed-complete}.  The valuations of $\pi$
at the six primes are
$(dN_0,dN_0,2N_0,2N_0(d-1),0,2dN_0)$.  The value $0$ occurs exactly
once, also for every power of $\pi$, so $\Q(\pi^m)=K'$ for all $m\ge1$
by \cref{lem:pair-valuations}.  All local degrees are one, so
\eqref{eq:HT-invariant-convention} gives the stated invariants, whose
least common denominator is $d$ because $d$ is even.  Honda--Tate theory
gives a division algebra of degree $d$ over $K'$ and a variety of
dimension $6d/2=3d$, and the base-change criterion gives geometric
simplicity and equality of the geometric and rational endomorphism
algebras.  The covariant slopes at the six primes are
$1/2,1/2,1-1/d,1/d,1,0$, each with multiplicity $d$, so the $p$-rank
is $d$.
\end{proof}

Choose a polarization of $B'_d$ defined over $\F_q$.  As in
\cref{sec:seeds}, its Rosati involution restricts to CM conjugation on
$K'$, because $\pi^\dagger=q/\pi=\bar\pi$ and $K'=\Q(\pi)$.  At the
first conjugate pair one has
$D'_{d,\mathfrak p_1}\simeq M_{d/2}(\Delta_2)$.

\begin{theorem}\label{thm:positive-rank}
Let $p$ be a prime, let $d\ge2$ be an even integer, and let $s>r\ge3$ be
coprime integers; put $n=r+s$.  There exist an algebraically closed
field $\Omega$ of characteristic $p$ and a simple abelian variety
$A/\Omega$ with
\[
 \dim A=3nd,\qquad f(A)=nd,\qquad m(A)=\frac nd,
\]
whose geometric endomorphism algebra is an algebra $D'_d$ as in
\cref{prop:sextic-seed}, for any sextic CM field $K'$ in which $p$
splits completely; such a field exists by \cref{lem:sextic-field}.
\end{theorem}

\begin{proof}
We follow the proof of \cref{thm:uniform-quartic}, replacing $B_d$ by
$B'_d$ and putting $A_0=(B'_d)^n$ with its product polarization.  Each
of the six local factors of $B'_d[p^\infty]$ has height $d$, so the
$\mathfrak p_1$-factor of $A_0[p^\infty]$ has height $nd$ and slope
$1/2$.  With the identification of $D'_{d,\mathfrak p_1}$ recorded
before the theorem, the Skolem--Noether argument of Step~1 of that
proof applies to $Y=X_2^{d/2}$.
\Cref{prop:quartic-globalization-complete}, now with $e_0=3$, gives a
polarized abelian scheme $\mathcal A/\overline{\F}_p[[t]]$ with rational
$D'_d$-action, special fibre $D'_d$-linearly $p$-power isogenous to
$A_0$, active factors $Y$ and $Y^\vee$, and constant factors at the
pairs $(\mathfrak p_2,\bar{\mathfrak p}_2)$ and
$(\mathfrak p_3,\bar{\mathfrak p}_3)$.  We check the hypotheses of
\cref{prop:commutant-criterion}.  The global ones hold by construction:
$D'_d$ acts rationally on $\mathcal A$, the special fibre is
$D'_d$-linearly isogenous to $(B'_d)^n$ with diagonal action, and
$\Endo(B'_d)=D'_d$ by \cref{prop:sextic-seed}.  The local ones concern
the $\mathfrak p_1$-factor $Y$ and are verified as in Step~2 of the
proof of \cref{thm:uniform-quartic}: the special fibre $Y_0$ is
isoclinic of slope $1/2$, the generic slopes $1/2-1/r$ and $1/2+1/s$
differ from $1/2$, and
\[
 \Endo_{D'_{d,\mathfrak p_1}}(Y_{\bar\eta})
 \simeq\mathcal D_r\times\mathcal D_s,\qquad \gcd(r,s)=1.
\]
Hence $\Endo(\mathcal A_{\bar\eta})=D'_d$, and $A=\mathcal A_{\bar\eta}$
is simple.

The active factors have slopes in $(0,1)$, and both slopes in the second
constant pair also lie in $(0,1)$.  The third pair is constant, with
slopes $1$ and $0$, each of height $nd$;
the slope-$1$ part is \'etale.  Hence $f(A)=nd$.  Finally, $\dim A=n\dim B'_d=3nd$ and
$m(A)=6nd/(6d^2)=n/d$.  Take $\Omega$ to be an algebraic closure of
$\overline{\F}_p((t))$.
\end{proof}

\begin{corollary}\label{cor:both-ranks}
For every prime $p$ and every $\mu\in\Q_{>0}$ there exist simple abelian
varieties $A$ and $A'$ over algebraically closed fields of
characteristic $p$ with $m(A)=m(A')=\mu$, $f(A)=0$ and $f(A')>0$.  In
particular, no positive rational number forces $p$-rank zero.
\end{corollary}

\begin{proof}
The variety $A$ is given by \cref{thm:main}.  For $A'$, if $\mu=1$, take
an ordinary elliptic curve over $\overline{\F}_p$; its endomorphism
algebra is an imaginary quadratic field, so $m=1$ and $f=1$.  If
$\mu=1/d$ with $d\ge2$, take $B^\circ_d\otimes\overline{\F}_p$ from
\cref{prop:unit-positive}.  Otherwise the reduced numerator $a$ of
$\mu=a/b$ satisfies $a\ge2$.  Choose an integer $c\ge1$ with $cb$ even
and $ca\ge7$, for instance $c=4$; put $n=ca$ and $d=cb$, and write
$n=r+s$ by \cref{lem:coprime-split}.  Since $d\ge2$ is even,
\cref{thm:positive-rank} gives $A'$ with $m(A')=\mu$, $\dim A'=3c^2ab$ and $f(A')=c^2ab$.
\end{proof}

\begin{corollary}\label{cor:r2-variant}
Let $d\ge4$ be even and let $s\ge3$ be odd.  For every prime $p$ there
exists a simple abelian variety $A$ over an algebraically closed field
of characteristic $p$ with
\[
 \Endo(A)=D_d,\qquad \dim A=2d(s+2),\qquad
 m(A)=\frac{s+2}d,\qquad f(A)=2d.
\]
In particular, if $a/b$ is a reduced fraction with $a$ odd and $b$
even, and $c$ is an odd integer with $ca\ge5$ and $cb\ge4$, then $d=cb$
and $s=ca-2$ realize $a/b$ in dimension $2c^2ab$ with $p$-rank $2cb$.
\end{corollary}

\begin{proof}
Take $r=2$ and $n=s+2$ in \cref{prop:two-level-window-complete}.  The
proof of \cref{thm:uniform-quartic} up to the determination of the
endomorphism algebra still applies: the active special slope is $1/2$,
the two generic slopes are $0$ and $1/2+1/s$, and the reduced
centralizer degrees $2$ and $s$ are coprime.  Both generic slopes differ
from the special slope, so \cref{prop:commutant-criterion} gives
$\Endo(A)=D_d$.  For $Y=X_2^{d/2}$, the slope-$0$ component has height
$2d$; its dual contributes slope $1$, that is, an \'etale part, with
multiplicity $2d$.  The other
slopes, including those of the constant pair, lie strictly between $0$
and $1$.  Thus $f(A)=2d$, and the dimension and ratio are computed as
before.  For the final assertion, $d=cb$ is even, $s=ca-2\ge3$ is odd,
and $(s+2)/d=a/b$.
\end{proof}

\begin{remark}\label{rem:comparison}
\Cref{thm:positive-rank} gives a uniform dimension bound, but smaller
realizations with positive $p$-rank are available in several ranges:
\cref{prop:unit-positive} for unit fractions (dimension $2d$, $p$-rank
$d$); \cref{prop:quadratic-realization} for reduced fractions $a/b$ with
$b+2\le a\le 2b-1$ (dimension $ab$, $p$-rank $b^2$);
\cref{cor:r2-variant} for reduced fractions $a/b$ with $a$ odd and $b$
even; and, for $\mu>2$, the Oort--van der Put construction of
\cref{cor:known-large-ratios}, whose $p$-rank is $2b^2$ in the notation
of its proof.  For $\mu=1/2$, a supersingular elliptic curve has
$p$-rank zero, while $B^\circ_2$ has dimension $4$ and $p$-rank $2$.
For $\mu=1$, an ordinary elliptic curve over $\overline{\F}_p$ gives
dimension $1$ and $p$-rank $1$. For $\mu=2$, an elliptic curve with transcendental $j$-invariant has
endomorphism ring $\Z$ and $p$-rank $1$.
For nonintegral ratios the ordinary case is excluded by
\cref{rem:ordinary}.
\end{remark}

\section*{Acknowledgements}
The author is deeply grateful to Professor Frans Oort for suggesting this problem, for his encouragement, and for his comments and questions on an earlier version of the manuscript, which led to improvements in the exposition. The author also thanks him for suggesting display theory and windows as tools for the deformation problem. 

This work was partially supported by CNPq (grant 302774/2025-4),
FAPESP (grant 2024/00923-6) and FAEPEX (grant 3485/25).

\bigskip
\noindent
\begin{minipage}{\textwidth}
\small
Institute of Mathematics, Statistics and Scientific Computing (IMECC),\\
University of Campinas (UNICAMP), Rua S\'ergio Buarque de Holanda 651,\\
13083-859 Campinas, SP, Brazil\\
\textit{Email address:} \href{mailto:saeed@unicamp.br}{\texttt{saeed@unicamp.br}}
\end{minipage}


\begin{thebibliography}{99}

\bibitem{BerthelotParfait}
P.~Berthelot,
\emph{Th\'eorie de Dieudonn\'e sur un anneau de valuation parfait},
Ann. Sci. \'Ec. Norm. Sup. (4) \textbf{13} (1980), no.~2, 225--268.

\bibitem{BBM}
P.~Berthelot, L.~Breen, and W.~Messing,
\emph{Th\'eorie de Dieudonn\'e cristalline~II},
Lecture Notes in Math., vol.~930, Springer-Verlag, Berlin, 1982.

\bibitem{ChaiMonodromy}
C.-L.~Chai,
\emph{Monodromy of Hecke-invariant subvarieties},
Pure Appl. Math. Q. \textbf{1} (2005), no.~2, 291--303.

\bibitem{ChaiOortMonodromy}
C.-L.~Chai and F.~Oort,
\emph{Monodromy and irreducibility of leaves},
Ann. of Math. (2) \textbf{173} (2011), no.~3, 1359--1396.

\bibitem{deJongHom}
A.~J.~de Jong,
\emph{Homomorphisms of Barsotti--Tate groups and crystals in positive
characteristic},
Invent. Math. \textbf{134} (1998), no.~2, 301--333;
erratum, Invent. Math. \textbf{138} (1999), no.~1, 225.

\bibitem{Demazure}
M.~Demazure,
\emph{Lectures on $p$-divisible groups},
Lecture Notes in Math., vol.~302, Springer-Verlag, Berlin, 1972.

\bibitem{FaltingsChai}
G.~Faltings and C.-L.~Chai,
\emph{Degeneration of Abelian Varieties},
Ergebnisse der Mathematik und ihrer Grenzgebiete (3), vol.~22,
Springer-Verlag, Berlin, 1990.

\bibitem{GilleSzamuely}
P.~Gille and T.~Szamuely,
\emph{Central Simple Algebras and Galois Cohomology},
Cambridge Studies in Advanced Mathematics, vol.~101,
Cambridge University Press, Cambridge, 2006.

\bibitem{EGAIII}
A.~Grothendieck,
\emph{\'El\'ements de g\'eom\'etrie alg\'ebrique: III. \'Etude cohomologique
 des faisceaux coh\'erents, Premi\`ere partie},
Publ. Math. Inst. Hautes \'Etudes Sci. \textbf{11} (1961), 5--167.

\bibitem{EGAIII2}
A.~Grothendieck,
\emph{\'El\'ements de g\'eom\'etrie alg\'ebrique: III. \'Etude cohomologique
 des faisceaux coh\'erents, Seconde partie},
Publ. Math. Inst. Hautes \'Etudes Sci. \textbf{17} (1963), 5--91.

\bibitem{EGAIV3}
A.~Grothendieck,
\emph{\'El\'ements de g\'eom\'etrie alg\'ebrique: IV. \'Etude locale
 des sch\'emas et des morphismes de sch\'emas, Troisi\`eme partie},
Publ. Math. Inst. Hautes \'Etudes Sci. \textbf{28} (1966), 5--255.

\bibitem{KatzST}
N.~M.~Katz,
\emph{Serre--Tate local moduli},
in \emph{Surfaces Alg\'ebriques (Orsay, 1976--78)},
Lecture Notes in Math., vol.~868, Springer-Verlag, Berlin, 1981, pp.~138--202.

\bibitem{Lau2014}
E.~Lau,
\emph{Relations between Dieudonn\'e displays and crystalline Dieudonn\'e theory},
Algebra Number Theory \textbf{8} (2014), no.~9, 2201--2262.

\bibitem{LenstraOort}
H.~W.~Lenstra Jr. and F.~Oort,
\emph{Simple abelian varieties having a prescribed formal isogeny type},
J. Pure Appl. Algebra \textbf{4} (1974), no.~1, 47--53.

\bibitem{Mumford}
D.~Mumford,
\emph{Abelian Varieties},
Tata Institute of Fundamental Research Studies in Mathematics, vol.~5,
Hindustan Book Agency, New Delhi, 2008; corrected reprint of the second
edition (1974).

\bibitem{GIT}
D.~Mumford, J.~Fogarty, and F.~Kirwan,
\emph{Geometric Invariant Theory}, third edition,
Ergebnisse der Mathematik und ihrer Grenzgebiete (2), vol.~34,
Springer-Verlag, Berlin, 1994.

\bibitem{OortEnd1988}
F.~Oort,
\emph{Endomorphism algebras of abelian varieties},
in \emph{Algebraic Geometry and Commutative Algebra in Honor of
Masayoshi Nagata}, vol.~II, Kinokuniya, Tokyo, 1988, pp.~469--502.

\bibitem{OortSomeQ}
F.~Oort,
\emph{Some questions in algebraic geometry},
June 1995, written for the conference ``Arithmetic and geometry of
abelian varieties''; reprinted as Appendix~1 of \cite{OortQ}.

\bibitem{OortNP}
F.~Oort,
\emph{Newton polygons and formal groups: conjectures by Manin and
Grothendieck},
Ann. of Math. (2) \textbf{152} (2000), no.~1, 183--206.

\bibitem{OortFF}
F.~Oort,
\emph{Abelian varieties over finite fields},
in \emph{Higher-Dimensional Geometry over Finite Fields}
(D.~Kaledin and Yu.~Tschinkel, eds.),
NATO Sci. Peace Secur. Ser. D Inf. Commun. Secur., vol.~16,
IOS Press, Amsterdam, 2008, pp.~123--188.

\bibitem{OortQ}
F.~Oort (ed.),
\emph{Open Problems in Arithmetic Algebraic Geometry},
Advanced Lectures in Mathematics, vol.~46, International Press, Somerville, MA,
2019;
Appendix~3, \emph{Questions in Arithmetic Algebraic Geometry}, pp.~295--331,
Questions~22.2 and~22.3 on p.~315.

\bibitem{OortVanDerPut}
F.~Oort and M.~van der Put,
\emph{A construction of an abelian variety with a given endomorphism algebra},
Compositio Math. \textbf{67} (1988), no.~1, 103--120.

\bibitem{Pierce}
R.~S.~Pierce,
\emph{Associative Algebras},
Graduate Texts in Mathematics, vol.~88, Springer-Verlag, New York, 1982.

\bibitem{Poonen}
B.~Poonen,
\emph{Rational Points on Varieties},
Graduate Studies in Mathematics, vol.~186,
American Mathematical Society, Providence, RI, 2017.

\bibitem{RZ}
M.~Rapoport and T.~Zink,
\emph{Period spaces for $p$-divisible groups},
Annals of Mathematics Studies, vol.~141, Princeton University Press, Princeton, NJ, 1996.

\bibitem{Tate66}
J.~Tate,
\emph{Endomorphisms of abelian varieties over finite fields},
Invent. Math. \textbf{2} (1966), 134--144.

\bibitem{WaterhouseMilne}
W.~C.~Waterhouse and J.~S.~Milne,
\emph{Abelian varieties over finite fields},
in \emph{1969 Number Theory Institute}, Proc. Sympos. Pure Math., vol.~20,
Amer. Math. Soc., Providence, RI, 1971, pp.~53--64.

\bibitem{ZinkSlope}
T.~Zink,
\emph{On the slope filtration},
Duke Math. J. \textbf{109} (2001), no.~1, 79--95.

\end{thebibliography}
\end{document}